\documentclass[a4paper,10pt]{amsart}
\usepackage{amsmath, amssymb, amsthm}
\numberwithin{equation}{section}
\newcommand{\hs}[1]{\hskip -#1pt}
\newcommand{\CR}{\mathbb{R}}
\newcommand{\CN}{\mathbb{N}}
\newcommand{\one}{1\!\!\!\;\mathrm{l}}
\newcommand{\eps}{\varepsilon}
\newcommand{\weak}{\rightharpoonup}
\newcommand{\half}{\frac{1}{2}}
\newcommand{\ahalf}{\frac{\alpha}{2}}
\newcommand{\loc}{\mathrm{loc}}

\newcommand{\T}{\mathcal{T}}
\newcommand{\A}{\mathcal{A}}
\newcommand{\expect}{\mathbb{E}}
\newcommand{\tr}{\mathrm{Tr}}
\newtheorem{thm}{Theorem}[section]
\newtheorem{prop}[thm]{Proposition}
\newtheorem{cor}[thm]{Corollary}
\newtheorem{lem}[thm]{Lemma}
\theoremstyle{definition}
\newtheorem{defn}[thm]{Definition}
\newtheorem{hyp}[thm]{Hypothesis}
\theoremstyle{remark}
\newtheorem{rem}[thm]{Remark}

\begin{document}

\title[Nonautonomous Kolmogorov parabolic equations]{Nonautonomous Kolmogorov parabolic equations with unbounded coefficients}

\author{Markus Kunze}
\address{Markus Kunze: Graduiertenkolleg 1100 Ulm University\newline
    Helmholtzstrasse 18, 89069 Ulm, Germany}
\email{markus.kunze@uni-ulm.de}

\author{Luca Lorenzi}
\address{L. Lorenzi, A. Lunardi: Dipartimento di Matematica, Universit\`a degli Studi di Parma, Viale G.P. Usberti, 53/A, I-43100 Parma, Italy}
\curraddr{}
\email{luca.lorenzi@unipr.it}
\email{alessandra.lunardi@unipr.it}
\thanks{Work supported by the M.I.U.R. research projects
Prin 2004 and 2006 ``Kolmogorov equations''.
To appear on Trans. Amer. Math. Soc.}

\author{Alessandra Lunardi}
\address{}
\curraddr{}
\email{}
\thanks{}


\begin{abstract}
We study a class of elliptic operators $A$ with unbounded coefficients defined in $I\times\CR^d$
for some unbounded  interval $I\subset\CR$. We prove that,
for any $s\in I$, the Cauchy problem $u(s,\cdot)=f\in C_b(\CR^d)$ for the parabolic equation $D_tu=Au$ admits a unique bounded
classical solution $u$.  This allows to
associate an evolution family $\{G(t,s)\}$ with $A$, in a natural way.
We study the main properties of this evolution family
and prove gradient estimates for the function $G(t,s)f$. Under suitable assumptions, we show that
there exists an evolution system of  measures for $\{G(t,s)\}$ and we study
the first properties of  the extension of $G(t,s)$ to
the $L^p$-spaces with respect to such measures.
\end{abstract}

\maketitle

\section{Introduction and summary}

Parabolic partial differential equations with unbounded coefficients occur naturally
in the study of stochastic processes. Let us consider the stochastic differential equation
\begin{equation}\label{sde}
\left\{\begin{array}{rcl}
dX_t &\hs{5} =\hs{5} & \mu (t, X_t)dt + \sigma (t, X_t) dW_t,\quad\;\,t>s,\\[1mm]
X_s &\hs{5} =\hs{5} & x.
\end{array}\right.
\end{equation}
Here, $W_t$ is a standard $d$-dimensional Brownian motion and $\mu$
(resp. $\sigma$) are regular $\CR^d$ (resp.
$\CR^{d\times d}$) valued coefficients. If (\ref{sde}) has a
solution $X_t = X(t,s,x)$ for all $x\in \CR^d$,  it follows from
It\^{o}'s formula that, for $\varphi \in C^2_b(\CR^d)$ and $t\in\CR$, the function
\[
u(s,x) := \expect \left( \varphi (X(t,s,x))\right)
\]
solves the partial differential equation
\begin{equation}\label{pde}
\left\{\begin{array}{rcl}
\displaystyle u_s(s,x) &\hs{5}=\hs{5} & - \half \tr ((\sigma(s,x) \sigma^*(s,x))
 D^2_xu(s,x)) - \langle \mu (s,x) , \nabla_x u(s,x) \rangle, \quad\;\,\, s < t,\\[2mm]
u(t,x) &\hs{5} =\hs{5} & \varphi (x).
\end{array}\right.
\end{equation}
This shows how probability theory may be used to obtain information
about the solutions of second order evolution PDE's.
In the case of Lipschitz continuous coefficients,
there are many results stating conditions on $\mu$ and $\sigma$ such
that (\ref{sde}) is well posed.
 See, e.g., \cite{ks,k,iw}.
\par
It is also possible to take (\ref{pde}) as a starting point and work in
a purely analytic manner.
This has been done in several papers in the autonomous case (see e.g., the book
\cite{lm} and its bibliography).
To the best of our knowledge, in the literature there is not any systematic treatment of the nonautonomous
case except in the particular case when the elliptic
operator in \eqref{pde} is the non autonomous Ornstein-Uhlenbeck operator
(see \cite{dpl,gl1,gl2}).

In this paper we set the basis for the general theory of non autonomous
operators. More precisely, we consider the equation
\begin{equation}\label{nonaueqn}
\left\{ \begin{array}{rcll}
\displaystyle  u_t(t,x) &\hs{5} =\hs{5} & \A (t)u(t,x), & (t,x) \in (s,+\infty)\times \CR^d,\\[2mm]
u(s,x) &\hs{5} =\hs{5} & f(x), & x\in \CR^d.
\end{array}\right.
\end{equation}
The operators $\A(t)$ appearing in \eqref{pde} are
defined on smooth functions $\varphi$ by
\begin{eqnarray*}
(\A(t)\varphi) (x) &\hs{5}=\hs{5}& \sum_{i,j=1}^d
q_{ij}(t,x)D_{ij}\varphi (x) + \sum_{i=1}^d b_i(t,x) D_i\varphi
(x)\\
&\hs{5} =\hs{5} & \tr \left( Q(t,x) D^2 \varphi (x)\right) + \langle b(t,x),
\nabla \varphi (x) \rangle.
\end{eqnarray*}
The time index $t$ varies over an interval $I$ which is either $\CR$ or a right halfline.

Note that the equation in \eqref{nonaueqn} is forward in time in contrast to equation (\ref{pde}). However,
reverting time, solutions of (\ref{nonaueqn}) are transformed into solutions
of (\ref{pde}) and viceversa.
Our standing hypotheses on the data $b = (b_i)$ and $Q = (q_{ij})$ are the following:
\begin{hyp}\label{hyp1}
\begin{itemize}
\item[(i)]
The coefficients $q_{ij}$ and $b_i$  belong to  $C^{\ahalf ,\alpha}_{\loc}(I\times \CR^d)$ for any
$i,j=1,\ldots,d$ and some $\alpha\in (0,1)$;
\item[(ii)] $Q$ is uniformly elliptic, i.e., for every $(t,x)\in I \times \CR^d$, the matrix $Q(t,x)$ is symmetric and there
exists a function $\eta:I\times\CR^d\to \CR$ such that $0<\eta_0:=\inf_{I\times\CR^d}\eta$ and
\[
\langle Q(t,x)\xi , \xi \rangle \geq \eta(t,x) |\xi |^2, \qquad\;\,
\xi \in \CR^d,\;\,(t,x) \in I \times \CR^d;
\]
\item[(iii)]
for every bounded interval $J \subset I$ there
exist a function $\varphi = \varphi_J \in C^2(\CR^d)$ and a positive
number $\lambda = \lambda_J$ such that
\[
\lim_{|x| \to +\infty}\varphi (x) = +\infty \quad \mbox{and}\quad
(\A(t)\varphi) (x) - \lambda \varphi (x) \leq 0,\qquad\;\,(t,x) \in
J\times \CR^d.
\]
\end{itemize}
\end{hyp}

Conditions (i) and (ii) are standard regularity and ellipticity assumptions in parabolic PDE's.
It is well known that assuming only (i) and (ii), problem \eqref{nonaueqn} may
admit several bounded solutions also in the autonomous case. Condition (iii)
is mainly used to ensure {\it uniqueness} of the bounded classical solution $u$ of \eqref{nonaueqn} (i.e.,
 uniqueness of a function $u\in C^{1,2}((s,+\infty)\times\CR^d)\cap C_b([s,T]\times\CR^d)$
 for any $T>s$, that satisfies \eqref{nonaueqn}).

\medskip
\par
In Section \ref{sect-2} we will be concerned with wellposedness of \eqref{nonaueqn} in the space
$C_b(\CR^d)$. In the autonomous case the solutions to \eqref{nonaueqn} are governed by a semigroup
$\{T(t)\}$ which is the transition semigroup of the Markov process obtained in (\ref{sde}).
In the non autonomous setting the semigroup is replaced by an evolution family $\{G(t,s)\}$.
We will establish several properties of this family in Section 3.
Note that, while regularity of $(G(t,s)\varphi)(x)$ with
respect to $(t,x)$ is a classical item in the theory of PDE's,
regularity with respect to $s$ is less standard. It is treated in the
literature in the case of bounded coefficients because of its
importance in several applications such as control theory. In our
case, to get continuity with respect to $s$ we have to sharpen
Hypothesis \ref{hyp1}(iii), assuming that $\A(t)\varphi$ is upperly
bounded in $J\times \CR^d$ for any bounded interval $J\subset I$.

\medskip
\par
In Section \ref{gradient-sect} we will study smoothing
properties of $G(t,s)$, proving several estimates on the spatial
derivatives of $G(t,s)\varphi$ for $\varphi\in C_{b}(\CR^d)$. We will
consider the following additional hypothesis:

\begin{hyp}\label{hyp2}
    \begin{itemize}
    \item[(i)]
    The data $q_{ij}$ and $b_i$ ($i,j=1,\ldots,d$) and their first-order spatial derivatives
belong to $C^{\ahalf , \alpha}_{\rm loc}(I\times\CR^d)$;
\item[(ii)]  there exists a continuous function $k : I \to
\CR$ such that
\[
\langle \nabla_x b(t,x) \xi, \xi \rangle \leq k(t) |\xi |^2,\qquad\;\,\xi\in\CR^d,\;\,(t,x)\in I\times\CR^d;
\]
\item[(iii)]
there exists a continuous function  $\rho : I \to [0, +\infty)$ such that,
for every $i,j,k \in \{1, \ldots , d\}$, we have
\[ |D_kq_{ij}(t,x)| \leq \rho (t) \eta (t,x),\qquad\;\,\xi\in\CR^d,\;\,(t,x)\in I\times\CR^d.
\]
\end{itemize}
\end{hyp}

Under this hypothesis we will prove uniform spatial gradient
estimates for the function $G(t,s)f$ when $f\in C^{k}_b(\CR^d)$, $k=0,1$, by means of
the classical Bernstein method (see \cite{B0}).
We will also prove more refined pointwise gradient estimates under either one of the following
more restrictive conditions:
\begin{hyp}
\label{hyp3}
\begin{itemize}
\item[(i)]
there exist a function $r:I\times\CR^d\to\CR$ and a constant $p_0\in (1,+\infty)$ such that
\[
\langle \nabla_x b(t,x) \xi, \xi \rangle \leq r(t,x) |\xi |^2,\qquad\;\,\xi\in\CR^d,\;\,(t,x)\in I\times\CR^d,
\]
and
$$
\sup_{(t,x)\in I\times\CR^d}\left(r(t,x) + \frac{d^3(\rho(t))^2\eta(t,x)}{4\min\{p_0-1,1\}}\right)<+\infty;
$$
\item[(ii)]
Hypothesis \ref{hyp2}(ii) holds true with the function $k$ being replaced by a real constant $k_0$. Moreover, there exists
a positive constant $\rho_0$ such that, for every $i,j,k=1,\ldots,d$, we have
\begin{eqnarray*}
|D_kq_{ij}(t,x)|\le\rho_0(\eta(t,x))^{\frac{1}{2}},\qquad\;\,(t,x)\in I\times\CR^d.
\end{eqnarray*}
\end{itemize}
\end{hyp}

Then, we get pointwise estimates,
\begin{equation}
|(\nabla G(t,s)f)(x)|^p \leq e^{\sigma_p(t-s)}(G(t,s)|\nabla f|^p)(x), \qquad\;\, t\geq s,\;\,x \in \CR^d,
\label{grad-punt0}
\end{equation}
for every $p\geq p_{0}$ and some real constant $\sigma_p$.
In the autonomous case (see \cite{lm}) these estimates are interesting not only for themselves,
but also for
the study of the behavior of the semigroup $\{T(t)\}$ in $L^{p}$-spaces
with respect to invariant measures.
An invariant measure corresponds to a stationary distribution of the Markov process
with transition semigroup $\{T(t)\}$. In the analytical setting, an
invariant measure for a Markov semigroup $\{T(t)\}$ is a Borel
probability measure such that
$$ \int_{\CR^d}(T(t)f)(y)\, \mu (dy) =  \int_{\CR^d}f(y)\, \mu (dy),
\quad t>0, \;f\in  C_b(\CR^d). $$
The interest in invariant
measures is due to the following:
\begin{enumerate}
\item[(i)]
 the invariant measure arises naturally in the asymptotic behaviour of the semigroup.
If $\mu$ is the (necessarily unique) invariant measure of $\{T(t)\}$, then
\[
(T(t)f)(x) \to \int_{\CR^d}f(y)\, \mu (dy)\quad \mbox{as}\,\,\, t \to +\infty,
\]
for any $f\in C_b(\CR^d)$ and any $x\in\CR^d$;
\item[(ii)]
the realizations of elliptic and parabolic operators in $L^p$-spaces with respect to invariant
measures are dissipative.
\end{enumerate}

In our nonautonomous case we cannot hope to find a
single invariant measure. Instead, we look for  systems of invariant
measures (see e.g. \cite{daprato,dynkin}) that is families of Borel
probability
measures $\{\mu_{t}: t\in I\}$ such that
$$ \int_{\CR^d}(G(t,s)f)(y)\, \mu_{t} (dy) =  \int_{\CR^d}f(y)\,
\mu_{s}(dy),
\quad t>s\in I, \;f\in  C_b(\CR^d). $$
In
Section \ref{evol}, we will prove the existence of a system of
invariant measures for our problem \eqref{nonaueqn} replacing Hypothesis \ref{hyp1}(iii)
with the following stronger condition:
\begin{hyp}\label{hyp4}
there exist a nonnegative function $\varphi \in
C^2(\CR^d)$, diverging to $+\infty$
as $|x|\to +\infty$, and constants $a, c > 0$ and $t_0 \in I$ such that
\[
(\A(t) \varphi) (x) \leq a - c\varphi (x),\qquad\;\,t\ge t_0,\quad\;\,x\in\CR^d.
\]
\end{hyp}
\medskip
\par
In contrast to the autonomous case, systems of invariant measures are, in general, not unique.
However, using a pointwise gradient estimate
we will prove that uniqueness holds in
the class of invariant measures $\{\mu_t: t\in I\}$ that admit finite moments
of  some order $p>0$,
which may blow up as $t\to +\infty$  with a certain
exponential rate. By definition, $\{\mu_t: t\in I\}$ admits finite moments
of order $p$ if, for any $t\in I$,
\begin{eqnarray*}
\mu_t(p)=\int_{\CR^d}|x|^p\mu_t(dx)<+\infty.
\end{eqnarray*}
Still using a uniform gradient estimate, we show that, also in the non autonomous case, the asymptotic behaviour
is determined by ``the'' system of invariant measures, in the sense that,
for any $x\in\CR^d$, $s\in I$,  and $ f\in  C_b(\CR^d)$,
\begin{eqnarray*}
\lim_{t\to +\infty} (G(t,s)f)(x) = \int_{\CR^d}f(y)\mu_s(dy),
\end{eqnarray*}
and the convergence is uniform in each compact set  in $\CR^d$.

Concerning point (ii), we note that, since we have to deal with
a family  of probability measures $\mu_{t}$,
we will also have a family of Lebesgue spaces $L^p(\mu_t)$
that are not mutually equivalent in general. This prevents us from extending
the operators $ G(t,s) $ to a single $L^p$-space, because $G(t,s)$
does not map $L^p(\mu_s)$ into itself in general, but it
maps $L^p(\mu_s)$ into $L^p(\mu_t)$. However,
it is possible to define an evolution semigroup associated
with $G(t,s)$ on a single $L^p$-space of functions defined in $I\times \CR^d$.
This was already done in \cite{dpl,gl1,gl2} in the special case
of Ornstein-Uhlenbeck operators.

The evolution semigroup associated with an evolution family $\{G(t,s)\}$ is
known to be a useful tool to determine several qualitative properties of the evolution
family. See e.g. the book \cite{ChiLa} and the references therein. In
the case of time depending Ornstein-Uhlenbeck operators, the use of
the evolution semigroup was essential to establish optimal regularity
results for evolution equations and also to get precise asymptotic
behavior estimates for $G(t,s)$, see \cite{gl1,gl2}.
However, the general theory of evolution semigroups is well established
only for evolution families acting on a fixed Banach space $X$, which is
not our case. Therefore, the study of the asymptotic behavior of
$G(t,s)\varphi $ for $\varphi \in L^p(\mu_s)$ through the evolution
semigroup is deferred to a future paper. Here, we just describe the first
properties of the evolution semigroup, in Section \ref{sect-5}.

In the last section we consider a simple example and see how our conditions may be
verified in this setting.

\subsection*{Notations}
We denote, respectively, by $ B_b(\CR^d)$ and $C_b(\CR^d)$ the set of all bounded
and Borel measurable functions $f:\CR^d\to\CR$ and its subset of all continuous functions.
We endow both   spaces with the sup-norm $\|\cdot\|_{\infty}$.

For any $k\in\CR_+$ (possibly $k=+\infty$) we denote by $C^k_b(\CR^d)$ the set of all
functions $f:\CR^d\to\CR$ that are continuously differentiable in $\CR^d$, up to $[k]$-th-order,
with bounded derivatives and such that the $[k]$-th-order derivatives are $(k-[k])$-H\"older continuous
in $\CR^d$. We norm $C^k_b(\CR^d)$ by
setting $\|f\|_{C^k_b(\CR^d)}=\sum_{|\alpha|\le [k]}\|D^{\alpha}f\|_{\infty}
+\sum_{|\alpha|=[k]}[D^{\alpha}f]_{C^{k-[k]}_b(\CR^d)}$.
$C^k_c(\CR^d)$ ($k\in\CN\cup\{+\infty\}$)  denotes the subset of $C^k_b(\CR^d)$ of
all compactly supported functions.  $C_0(\CR^d)$ denotes the set of all continuous functions
vanishing at infinity.

If $f$ is smooth enough we set $D_jf=\frac{\partial f}{\partial x_j}$,
\begin{eqnarray*}
|\nabla f(x)|^2=\sum_{i=1}^d|D_if(x)|^2,\qquad\;\,
|D^2f(x)|^2=\sum_{i,j=1}^d|D_{ij}f(x)|^2
\end{eqnarray*}
and
\begin{eqnarray*}
\|\nabla f\|^2_{\infty}=\sup_{x\in\CR^d}|\nabla f(x)|^2,\qquad\;\,
\|D^2f\|^2_{\infty}=\sup_{x\in\CR^d}|D^2 f(x)|^2.
\end{eqnarray*}

Suppose that $f$ depends on both time and spatial variables. If there is damage of confusion,
 we denote by $\nabla_xf$ and $D^2_xf$
the gradient and the Hessian matrix of the function $f(t,\cdot)$.
When $f$ is a vector valued function, $\nabla_xf$ denotes the Jacobian matrix of $f(t,\cdot)$.

Let $D\subset\CR^{d+1}$ be a domain or the closure of a domain.
By $C^{k+\alpha/2,2k+\alpha}_{\rm loc}(D)$
($k=0,1$, $\alpha\in (0,1)$)
we denote the set of all functions $f:D\to\CR$ such that the time derivatives
up to the $k$-th-order
and the spatial derivatives up to the $2k$-th-order are H\"older continuous with
exponent $\alpha$, with respect to the parabolic distance, in any compact set $D_0\subset D$.

For any $r>0$ we denote by $B_r \subset\CR^d$ the open ball centered at
$0$ with radius $r$.

Given a measurable set $E$, we denote by \mbox{$\one_E$} the characteristic function of $E$,
i.e., \mbox{$\one_E(x)=1$} if $x\in E$, $\one_E(x)=0$ otherwise.

Finally, we use the notation  $u_f$ for the (unique) bounded classical solution to  problem \eqref{nonaueqn}.


\section{Solutions in $C_b(\CR^d)$}
\label{sect-2}

In this section we want to solve our parabolic problem \eqref{nonaueqn} with data
$s \in I$ and $f \in C_b(\CR^d)$. By a {\it solution} of \eqref{nonaueqn} we mean a bounded classical
solution, i.e. a
function $u \in C_b([s, +\infty) \times
\CR^d) \cap C^{1 , 2 }( (s , +\infty ) \times
\CR^d)$ such that \eqref{nonaueqn} is satisfied. In the whole section
we
assume that Hypotheses \ref{hyp1} is fulfilled.

We already mentioned that
Hypothesis \ref{hyp1}(iii)  ensures uniqueness of the solution to \eqref{nonaueqn}.
In fact, it implies a maximum principle that we state as our first

\begin{thm}\label{maximum}
Let $s\in I$ and $T > s$. If $u \in C_b([s,T]\times
\CR^d)\cap C^{1,2}( (s,T] \times \CR^d)$ satisfies
\[ \left\{\begin{array}{rclr}
u_t(t,x) - \A(t)u(t,x) &\hs{5} \leq\hs{5} &  0, & \,\, (t,x)\in
(s,T]\times\CR^d,\\[1.5mm]
u(s,x) &\hs{5} \leq\hs{5} & 0, & x \in \CR^d,
\end{array}\right.
\]
then $u \leq 0$.
\end{thm}

\begin{proof}
The proof can be obtained as the proof of \cite[Theorem 4.1.5]{lm}.
\end{proof}

\par
Now, let us prove that problem \eqref{nonaueqn} admits a unique bounded
classical solution for any $f\in C_b(\CR^d)$.

\begin{thm}\label{thm1}
For every $s \in I$ and
every $f \in C_b(\CR^d)$, there exists a unique solution $u$ of problem
\eqref{nonaueqn}. Furthermore,
\begin{equation}
\|u(t,\cdot)\|_{\infty} \leq
\|f\|_{\infty},\qquad\;\, t \geq s.
\label{contractive}
\end{equation}
\end{thm}

\begin{proof}
Uniqueness follows from applying Theorem \ref{maximum} to $u-v$ and to
$v-u$, if $u$ and $v$ are two solutions.
Estimate \eqref{contractive} follows by applying
the same theorem to $\pm  u- \|f\|_{\infty} $.

The existence part can be obtained in a classical way, solving
Cauchy-Dirichlet problems in the balls $B_{n}$ and then letting $n\to
+\infty$. See e.g., \cite[Proposition 2.2]{simona-giorgio-enrico},
 \cite[Theorem 4.2]{mpw}, \cite[Theorems 2.2.1, 11.2.1]{lm}.
Since there are some technicalities, for the reader's convenience we go into details.
We split the proof in three steps.
\medskip
\par
{\em Step 1.} Here, we assume that $f$ belongs to $C^{2+\alpha}_c(\CR^d)$.
Denote by $n_0$ the smallest integer such that $\mathrm{supp}(f)$ is contained in the ball
$B_{n_0}$. Further, for any $n\ge n_0$, we consider the Cauchy-Dirichlet problem
\begin{equation}
\left\{ \begin{array}{rclll}
u_t(t,x) &\hs{5} =\hs{5} & \A(t)u(t,x), & t\in (s, +\infty ), &x\in B_n,\\[1mm]
u(t,x) &\hs{5} =\hs{5} & 0, & t\in (s,+\infty), &x \in
\partial B_n,\\[1mm]
u(s,x) &\hs{5} =\hs{5} & f(x), && x \in B_n, \end{array}\right.
\label{approx-dirichlet}
\end{equation}
in the ball $B_n$.
By classical results (see e.g., \cite{friedman} or \cite{lsu}) and Hypotheses \ref{hyp1}(i)-(ii), for any $n\ge n_0$,
the problem \eqref{approx-dirichlet} admits a unique solution
$u_n\in C^{1+\ahalf , 2+
\alpha}_{\loc}([s,+\infty)\times \overline{B_n})$.
Moreover, the classical Schauder estimates imply that,
for any $m\in\mathbb N$, with $m>n_0$, there exists a
constant $C=C(m)$ independent of $n$, such that
\begin{equation}
\|u_n\|_{C^{1+\ahalf ,2+\alpha}(D_m)} \leq C\|f\|_{C^{2+\alpha}(\CR^d)},
\label{uptobound}
\end{equation}
for any $n>m$, where $D_m=(s,m)\times B_{m}$. By the Arzel\`{a}-Ascoli
theorem, there exists a subsequence $(u^m_n)$ of $u_n$ which converges in $C^{1,2}(\overline{D_m})$ to some function $u^m\in
C^{1+\ahalf , 2+\alpha}(D_m)$. Of course, $u^m$ satisfies the differential equation $D_t u^m = {\mathcal A}(\cdot)u^m$ in $D_m$
and it equals $f$ on $\{s\}\times B_{m}$. Without loss of generality, we can assume that $u^{m+1}_n$
is a subsequence of $u^m_n$.
Note that, in this case, $u^{m+1}|_{D_m} \equiv u^m$. Hence, we can
define a function $u$ by putting $u|_{D_m} := u^m$. A standard procedure shows that $u$ belongs
to $C^{1+\alpha/2,2+\alpha}_{\rm loc}([s,+\infty)\times\CR^d)$ and it is the classical solution
to problem \eqref{nonaueqn}.

Note that the sequence $u_n$ itself converges to $u$ as $n$ tends to $+\infty$,
locally uniformly in $[s,+\infty)\times\CR^d$. Indeed, the above arguments
show that any convergent subsequence of $(u_n)$ should converge to a classical solution of \eqref{nonaueqn}.
\medskip
\par
{\em Step 2.}
Assume now that $f\in C_0(\CR^d)$. Then, there exists a sequence $(f_n)\subset C^{2+\alpha}_c(\CR^d)$ converging
to $f$ uniformly in $\CR^d$ as $n$ tends to $+\infty$.
Estimate \eqref{contractive} yields
\begin{eqnarray*}
\|u_{f_n}-u_{f_m}\|_{C_b([s,+\infty)\times\CR^d)}\le \|f_n-f_m\|_{C_b(\CR^d)},\qquad\;\,m,n\in\CN.
\end{eqnarray*}
Therefore, there exists a bounded and continuous function $u$ such that $u_{f_n}$ converges
to $u$, uniformly in $[s,+\infty)\times\CR^d$.
Moreover, applying the interior Schauder estimates to the sequence $(u_{f_n})$, we deduce that
$u_{f_n}$ converges in $C^{1,2}_{\rm loc}((s,+\infty)\times\CR^d)$ to $u$.
Hence, $u$ is the bounded classical solution of
problem \eqref{nonaueqn}.
\medskip
\par
{\em Step 3.} Now, fix $f\in C_b(\CR^d)$ and consider a bounded sequence
$(f_n)\in C^{2+\alpha}_c(\CR^d)$ converging to $f$ locally uniformly in $\CR^d$
as $n$ tends to $+\infty$. The same arguments as in Step 2 show that, up to a subsequence, $u_{f_n}$ converges,
in $C^{1,2}_{\rm loc}((s,+\infty)$ $\times\CR^d)$,
to some function $u\in C^{1+\alpha/2,2+\alpha}_{\rm loc}((s,+\infty)\times\CR^d)$,
as $n$ tends to $+\infty$. In particular, $u$ solves the differential equation in \eqref{nonaueqn}.
To prove that $u$ is, actually, a
classical solution of the problem \eqref{nonaueqn}, we fix a compact set
$K\subset\CR^d$ and a smooth and compactly supported function $\varphi$
such that $0\le\varphi\le 1$ and $\varphi\equiv 1$ in $K$. Further, we split $u_{f_n}=u_{\varphi f_n}+
u_{(1-\varphi)f_n}$, for any $n\in \CR$.
Since the function $\varphi f$ is compactly supported in $\CR^d$, it follows from Step 2  that
$u_{\varphi f_n}$ converges to $u_{\varphi f}$ uniformly in $[s,+\infty)\times\CR^d$.

Let us now consider the sequence $(u_{(1-\varphi)f_n})$. Fix $m \in \CN$. We claim that
\begin{equation}
\label{est1}
|(u_{(1-\varphi ) f_m})(t,x)| \leq (1 - u_{\varphi}(t,x))M,\qquad\;\,(t,x)\in (s,+\infty)\times\CR^d,
\end{equation}
where $M=\sup_{n\in\CN}\|f_n\|_{\infty}$.
Indeed, as a straightforward computation shows,
\[
(1 - u_{\varphi})M=u_{(1-\varphi)M}.
\]
Therefore, the function $w:=u_{(1-\varphi)f_m} -M(1-u_{\varphi})$
satisfies $w_t={\mathcal A}w$ and, moreover,
$$
w(s,\cdot)=(1-\varphi)f_m-M(1- \varphi )=(1-\varphi)(f_m-M)\le 0.
$$
The maximum principle of Theorem \ref{maximum} immediately implies that
$w$ is nonpositive in $[s,+\infty) \times\CR^d$ or, equivalently, $u_{(1-\varphi ) f_m}\le
(1 - u_{\varphi})M$. To prove the other inequality in \eqref{est1}, it suffices to observe that
$(- u_{(1-\varphi ) f_m}) =(u_{-(1-\varphi ) f_m})$ and repeat the above arguments with $f_m$
 replaced by $-f_m$.

Now, since $u_{f_n}$ converges pointwise to $u$, for any $(t,x)\in (s,+\infty)\times\CR^d$ we have
\[
|u(t,x)-f(x)| =\lim _{n\to +\infty} |u_{ f_n}(t,x)-f( x)|,\qquad\;\,(t,x)\in (s,+\infty)\times\CR^d,
\]
and, for each $n\in \CN$, we have
\begin{eqnarray}
    |u_{ f_n}(t,x)-f(x)| &\hs{5}\le\hs{5}&
|u_{ \varphi f_n}(t,x) -f(x)| + |u_{(1-\varphi)f_n}(t,x)|\nonumber\\
&\hs{5}\leq\hs{5}& |u_{\varphi f_{n}}(t,x)-f(x)|+
(1 - u_{\varphi}(t,x))M.
\label{continuity}
\end{eqnarray}
The right-hand side of
\eqref{continuity} converges to $0$ uniformly in $K$ as $t$ tends to $s^+$. Hence, $u$ can be
continuously extended  up to $t=s$   setting $u(s,\cdot)=f$. This completes the proof.
\end{proof}

\begin{rem}\label{rem2.3}
Let us observe that the choice of approximating   problem \eqref{nonaueqn}
by  Cauchy-Dirichlet problems in the ball $B_n$ is not essential. Indeed, repeating step by
step the proof of Theorem \ref{thm1}, we can see that   problem \eqref{nonaueqn} can be
approximated also by the Cauchy-Neumann problems
\begin{equation}
\left\{ \begin{array}{rcll}
u_t(t,x) &\hs{5} =\hs{5} & \A(t)u(t,x), & (t,x) \in (s,+\infty )\times B_n,\\[2mm]
\displaystyle\frac{\partial}{\partial\nu}u(t,x) &\hs{5} =\hs{5} & 0, & (t,x)\in (s,+\infty)\times\partial B_n,\\[2mm]
u(s,x) & \hs{5}= \hs{5} & f(x), &x\in B_n.
\end{array}\right.
\label{approx-neumann}
\end{equation}
We will use this approach in Section \ref{gradient-sect} to prove  estimates
for the space derivatives of $G(t,s)f$.
\end{rem}

Now we define the {\it evolution family} associated with our problem
\eqref{nonaueqn}. Let
\begin{eqnarray*}
\Lambda := \{(t,s) \in I\times I \, : \, t \geq s\}.
\end{eqnarray*}
 We put
$G(t,t) := id_{C_b(\CR^d)}$ and for $t > s$ we define the operator
$G(t,s)$ by setting
\[
(G(t,s)f)(x) := u_{f}(t,x),\qquad\;\, x\in\CR^d,
\]
where $u_f$ is the unique solution of problem \eqref{nonaueqn}. We call the
family $\{G(t,s): (t,s)\in\Lambda \}$ the {\it evolution family
associated with the problem} \eqref{nonaueqn}.

It is immediate from Theorem \ref{maximum} that, for $(t,s) \in
\Lambda$, the operator $G(t,s)$ is a positive contraction on
$C_b(\CR^d)$. From the uniqueness assertion in Theorem \ref{thm1},
the {\it law of evolution}
\begin{equation}\label{loe}
G(t,s)G(s,r) = G(t,r),
\end{equation}
for $r\leq s \leq t$ easily follows.

\hs{1}
The connection with Markov processes suggests that every operator
$G(t,s)$ should be associated with a transition kernel. Recall that
a transition kernel $p$ is a mapping from $\CR^d\times \mathcal{B}(\CR^d)$
such that $p(x, \cdot )$ is a sub probability measure for fixed $x$ and
such that $p(\cdot , A)$ is measurable for fixed $A\in \mathcal{B}(\CR^d)$.
The following proposition states that this is indeed the case; in fact, the transition
kernels $p_{t,s}$ form the non autonomous
equivalent of a conservative, stochastically continuous transition function,
cf. \cite[Sections 2.1 and 2.8]{dynkin}.

\begin{prop}\label{transition}
For every $(t,s)
\in \Lambda$ and every $x \in \CR^d$ there exists a unique
probability measure $p_{t,s}(x, \cdot )$ such that
\begin{equation}\label{intrep}
(G(t,s)f)(x) = \int_{\CR^d} f(y) \, p_{t,s}(x,dy),
\end{equation}
for each $C_b(\CR^d)$.
Furthermore, the following properties hold:
\begin{enumerate}
\item[(i)] for every $t \in I$, $p_{t,t}(x, \cdot)$ is
the Dirac measure concentrated at $x$;
\item[(ii)] for $t>s$ the measure $p_{t,s}(x, \cdot )$ is equivalent to the
Lebesgue measure, i.e. they have the same sets of zero measure;
\item[(iii)] for fixed $A \in B_b(\CR^d)$ and $(t,s) \in \Lambda$
the map $x \mapsto p_{t,s}(x,A)$ is Borel measurable;
\item[(iv)] for every open set $U \subset \CR^d$ containing $x$ we have
\[
\lim_{t\to s^+}p_{t,s}(x,U) = 1;
\]
\item[(v)] for every $t \geq s \geq r \, ,\, x \in \CR^d$ and $A \in
\mathcal{B}(\CR^d)$ we have
\[p_{t,r}(x,A) = \int_{\CR^d} p_{s,r}(y,A) \, p_{t,s}(x,dy).
\]
\end{enumerate}
\end{prop}

\begin{proof}
Let us define
\[
p_{t,s}(x, A) = \int_A g(t,s,x,y) \, dy,
\]
for any (Lebesgue) measurable set $A\subset\CR^d$ and any $t>s$, where
$g$ is the Green function of problem \eqref{nonaueqn} which can be obtained as the pointwise limit of the increasing
(with respect to $n$) sequence of Green functions $g_n$ associated with
problem \eqref{approx-dirichlet}. For the existence of these latter kernels,
see e.g., \cite[Theorem 3.16]{friedman}.
The function $g$ is measurable in its entries and it is positive since the $g_n$'s are.

Notice that, since \mbox{$G(t,s)\one \equiv \one$} by uniqueness, we have
\[
p_{t,s}(x,\CR^d) = \int_{\CR^d}\one (y)\,p_{t,s}(x,dy) =
(G(t,s)\one) (x) = 1,
\]
i.e. $p_{t,s}(x,\cdot)$ is a probability
measure.
Now formula  \eqref{intrep} and properties (i) and (ii) immediately follow.

To prove (iii), let $A \in
\mathcal{B}(\CR^d)$. Then, there exists a bounded sequence $(f_n)\subset
C_b(\CR^d)$ converging almost everywhere to $\one_A$. Hence, by the
dominated convergence theorem,
\begin{eqnarray*}
p_{t,s}(x, A) &\hs{5}=\hs{5}& \int_{\CR^d}\one_A(y) \, p_{t,s}(x,dy)\\
&\hs{5}=\hs{5}& \lim_{n\to +\infty} \int_{\CR^d} f_n(y) \, p_{t,s}(x, dy)\\
&\hs{5} =\hs{5}&
\lim_{n\to +\infty} (G(t,s)f_n)(x),
\end{eqnarray*}
for any $(s,t)\in\Lambda$ and any $x\in\CR^d$, and this implies that the function $(t,s,x) \mapsto p_{t,s}(x,A)$ is measurable.

Property (iv) follows from the continuity of the map $G(\cdot,s)f$ on $\{s\}\times\CR^d$, by virtue of
\cite[Lemma 2.2]{dynkin}.

Finally, (v) is an immediate consequence of (\ref{loe}).
\end{proof}

\begin{cor}
\label{irreducibility}
The evolution family $\{G(t,s)\}$ is irreducible, i.e.,
$(G(t,s)\one_{U})(x)>0$ for any open set $U\subset\CR^d$, any $(s,t)\in\Lambda$ and any $x\in\CR^d$.
More generally, \mbox{$(G(t,s)\one_A)(x) >0$} for each Borel set $A\subset \CR^d$
with positive Lebesgue measure, and for any $(s,t)\in\Lambda$ and any $x\in\CR^d$.
\end{cor}

\begin{proof} The   statement follows from the inequality $g(t,s,x,y)>0$ ($(s,t)\in\Lambda$, $x\in\CR^d$)
which holds since the Green function $g$ is the pointwise limit  of the increasing
  sequence of Green functions $g_n$ associated with
problem \eqref{approx-dirichlet}.
\end{proof}
\begin{rem}
The representation (\ref{intrep}) implies that we can extend our evolution family $G(t,s)$
to the space $B_b(\CR^d)$ of bounded, measurable functions.
More generally, in the sequel we set
\[
(G(t,s)f)(x) := \int_{\CR^d}f(y) p_{t,s}(x,dy),\qquad\;\,x\in\CR^d,
\]
also for unbounded functions $f:\CR^d\to\CR$, provided
that $f\in L^1(\CR^d, p_{t,s}(x, \cdot))$.

Formula
(\ref{intrep}) also implies that
the adjoints $G^*(t,s)$ leave the space of signed measures invariant.
\end{rem}


\section{Continuity properties of the evolution family $\{G(t,s)\}$}

In this section we prove some useful continuity properties of the function $G(t,s)f$ when
$f\in C_b(\CR^d)$.
To begin with, let us prove the following proposition.

\begin{prop}\label{properties}
Let $(f_n) \subset
C_b(\CR^d)$ be a   sequence such that $\|f_n\|_{\infty} \leq M$ for each $n\in \CN$. Then:
\begin{enumerate}
\item[(i)]
if $f_n$ converges pointwise to $f$, then $G(\cdot , s )f_n$ converges
to $G(\cdot , s )f$ uniformly in $[s+\eps,T]\times K$ for each $T>s+\eps>s$ and each compact
set $K \subset \CR^d$;
\item[(ii)]
 if $f_n$
converges uniformly to $f$ in each compact set $K\subset\CR^d$, then $G(\cdot,s)f_n$
converges to $G(\cdot , s )f$ uniformly in
$[s,T] \times K$, for each $T>s$ and each compact set $K \subset
\CR^d$.
\end{enumerate}
\end{prop}

\begin{proof}
(i). By the representation formula (\ref{intrep}) and
the dominated convergence theorem, we obtain that $G(\cdot,s)f_n$ converges to $G(\cdot,s)f$
pointwise in $[s,+\infty)\times\CR^d$, as $n$ tends to $+\infty$. However, the
classical interior Schauder estimates yield
\begin{eqnarray*}
\|G(\cdot,s)f_n\|_{C^{1+\alpha/2,2+\alpha}([s+\varepsilon,T]\times K)}\le
C \|f_n\|_{\infty}\le CM, \qquad\;\,n\in\CN,
\end{eqnarray*}
for any $\varepsilon>0$, any compact set $K\subset\CR^d$ and
some positive constant $C=C(\varepsilon,K)$.  The Arzel\`a-Ascoli theorem implies that
$G(\cdot,s)f_n$ converges to $G(\cdot,s)f$ in $C^{1,2}([s+\varepsilon,T]\times K)$.
\medskip
\par
(ii). The proof follows adapting the arguments used in the proof of
Theorem \ref{thm1}.
Let $K$ be a compact set and let
$\varphi\in C^{2+\alpha}_c(\CR^d)$ be such that $\varphi\equiv 1$ in $K$.
Split $u_{f_n}=u_{\varphi f_n}+u_{(1-\varphi)f_n}$. By Step 2 in the proof of Theorem \ref{thm1},
$u_{\varphi f_n}$ converges to $u_{\varphi f}$ uniformly in $[s,+\infty)\times\CR^d$.

To complete the proof, it suffices to show that $u_{(1-\varphi)f_n}$ converges
to $u_{(1-\varphi)f}$ uniformly in $[s,T]\times K$ as $n$ tends to $+\infty$ for any $T>s$.
The arguments in Step 2 of the proof of Theorem \ref{thm1} show
that, up to a subsequence, $u_{(1-\varphi)f_n}$ converges in $C^{1,2}([s+\varepsilon,T]\times B_r)$,
for any $\varepsilon>0$ and any $r>0$,
to a function $v\in C^{1+\alpha/2,2+\alpha}_{\rm loc}((s,+\infty)\times\CR^d)$.
Moreover, letting $m$ go to $+\infty$ in \eqref{est1} gives
\begin{equation}
|v(t,x)| \leq (1 - u_{\varphi}(t,x))M,\qquad\;\,(t,x)\in (s,+\infty)\times\CR^d.
\label{stima-utile}
\end{equation}
Since $u_{\varphi}$ is continuous at $\{s\}\times\CR^d$ and $\varphi\equiv 1$ in $K$, from
\eqref{stima-utile}  we   deduce that
$v(t,x)$ converges to $0$ as $t\to s^+$, uniformly with respect to $x\in K$.

Let us now fix $\eps>0$ and let $\delta$ be sufficiently small such that
$|u_{(1-\varphi)f_n}|+|v|\le\eps/2$ in $[s,s+\delta]\times K$ for any $n\in\CN$.
Moreover, we fix $n$ large enough such that $|u_{(1-\varphi)f_n}-v|\le\eps/2$ in $[s+\delta,T]\times K$.
For such $n$ and $\delta$ we get
\begin{eqnarray*}
\|u_{(1-\varphi)f_n}-v\|_{C([s,T]\times K)}
&\hs{5} \leq\hs{5} & \|u_{(1-\varphi)f_n}\|_{C([s,s+\delta]\times K)} +\|v\|_{C([s,s+\delta]\times K)}\nonumber\\
&\hs{5}&+\|u_{(1-\varphi)f_n} - v\|_{C[s+\delta , T]\times
K)}\\
&\hs{5} \leq\hs{5} & \eps.
\end{eqnarray*}
Summing up,  the sequence  $u_{f_n}=u_{\varphi f_n}+u_{(1-\varphi)f_n}$ converges
as $n$ tends to $+\infty$,
to the function $u=u_{\varphi f}+v$ which belongs to
$C^{1+\alpha/2,2+\alpha}_{\rm loc}((s,+\infty)\times\CR^d)$, and the convergence is uniform in
$[s,T]\times K$. Since  $K$ is arbitrary, $u_{f_n}$ converges locally uniformly to
$u$ in $[s,+\infty)\times\CR^d$, so that  $u$ is continuous up to $t=s$ where it equals $f$.
Moreover, since $u_{f_n}$
converges to $u$ in $C^{1,2}([s+\varepsilon,T]\times B_R)$ for any
$\varepsilon\in (0,T-s)$ and any $R>0$,
then $D_tu-{\mathcal A}u=0$ for $t>s$.
Thus, $u$ is a bounded classical solution of \eqref{nonaueqn}
and, by Theorem \ref{maximum}, $u=u_f$.
This completes the proof.
\end{proof}

\subsection{Continuity of the function $G(t,s)f$ with respect to the variable $s$}

Since evolution families depend on two parameters $t$ and $s$, it is natural to investigate
 also the smoothness of the function $G(t,\cdot)f$. In the
following lemma we prove a very useful generalization of the well known
formula that holds in the   case of bounded coefficients.
This lemma will play a fundamental role to prove existence of evolution systems of invariant measures
in Section \ref{evol}.

\begin{lem}\label{sderivative}
Let $f \in C_{b}^2(\CR^d)$ be constant outside a compact set $K$.
Then, for any $x \in \CR^d$, and any $s_0 <s_1 \leq  t$, the function $ r
\mapsto (G(t,r)\A(r)f)(x)$ is integrable in $(s_{0},s_{1})$ and we
have
\begin{equation}
(G(t,s_1)f)(x) - (G(t,s_0)f)(x) = - \int_{s_0}^{s_1}
(G(t,r)\A(r)f)(x) \, dr.
\label{dom-int}
\end{equation}
In particular,
the function $(G(t,\cdot)f)(x)$ is continuously differentiable in
$I_t:=I\cap (-\infty,t]$ and
\begin{equation}
\frac{d}{ds}(G(t,s)f)(x) = -(G(t,s)\A(s)f)(x),\qquad\;\,s\in I_t.
\label{deriv}
\end{equation}
Finally for any $g \in C_0(\CR^d)$ the function $G(t,\cdot )g$ is continuous in $I_t$ with
values in $C_b(\CR^d)$.
\end{lem}

\begin{proof}
By assumption, we can
write $f= g +c\cdot \one$, for some $g\in C^2_c(\CR^d)$ and some $c\in\CR$. However,
\mbox{$G(t,s)\one \equiv \one$}, whence the assertion is trivially satisfied by any
constant function.
Thus, it remains to prove it when $f\in C_c^2(\CR^d)$.
Choose $n_0$ such that $\mathrm{supp}(f) \subset B_{n_0}$, and denote by
$\{G_n(t,s)\}$ the evolution family associated with problem \eqref{approx-dirichlet} for $n \geq
n_0$ (cf. \cite[Theorem 6.3]{acqui}).
By \cite[Theorem 2.3(ix)]{acqui}, we can write \eqref{deriv} with $G$ being replaced by $G_n$.
Integrating such an equality with respect to $s$ and  recalling that, by Step 1 in the proof of Theorem \ref{thm1},
for any $(t,r)\in\Lambda$, $G_n(t,r)f$ converges to $G(t,r)f$ pointwise in $\CR^d$ as $n$ tends to $+\infty$,
we obtain
\begin{eqnarray}
(G(t,s_1)f)(x) - (G(t,s_0)f)(x) &\hs{5} =\hs{5} & \lim_{n\to +\infty} (G_n(t,s_1)f)(x) -
(G_n(t,s_0)f)(x)\nonumber \\
&\hs{5} =\hs{5} & - \lim_{n\to +\infty}  \int_{s_0}^{s_1} (G_n(t,r)\A(r)f)(x) \, dr\nonumber\\
&\hs{5} =\hs{5} & - \int_{s_0}^{s_1} (G(t,r)\A(r)f)(x)\, dr,
\label{cont}
\end{eqnarray}
where the last equality follows by dominated convergence.

Now, observe that
\eqref{cont} implies that the function $G(t,\cdot)f$ is continuous in $I_t$,
with values in $C_b(\CR^d)$, for any $f\in C^2_c(\CR^d)$. Since $C^2_c(\CR^d)$ is dense in $C_0(\CR^d)$,
$G(t,\cdot)g$ is continuous in $I_t$ for any $f\in C_0(\CR^d)$.

To prove that the function $G(t,\cdot)f$ is differentiable,
it is enough to show that the function $G(t,\cdot)\A(\cdot)f$ is continuous in
$I_t$. Indeed, for any $r,r_0\in I_t$,
\begin{eqnarray*}
&\hs{5}&\|G(t,r)\A(r)f-G(t,r_0)\A(r_0)f\|_{\infty}\nonumber\\
&\hs{5}=\hs{5}&\|G(t,r)(\A(r)f-\A(r_0)f)\|_{\infty}+\|(G(t,r)-G(t,r_0))\A(r_0)f\|_{\infty}\nonumber\\
&\hs{5}\le\hs{5}&\|\A(r)f-\A(r_0)f\|_{\infty}+\|(G(t,r)-G(t,r_0))\A(r_0)f\|_{\infty},
\end{eqnarray*}
and the last side of the previous chain of inequalities goes to $0$ as $r\to r_0$, since
$\A(r_0)f\in C_c(\CR^d)$.
Now, \eqref{cont} implies that the function
$G(t,\cdot)f$ is differentiable in $I_t$, and \eqref{deriv} follows.
This completes the proof.
\end{proof}

To prove that $(t,s,x)\mapsto (G(t,s)f)(x)$ is continuous in $\Lambda\times\CR^d$
for any function
$f\in C_b(\CR^d)$, we need
an intermediate assumption between Hypothesis \ref{hyp1}(iii) and Hypothesis \ref{hyp4}.
More precisely, in the rest of this section we assume that the following hypothesis is satisfied.

\begin{hyp}\label{hyp5}
For every bounded interval $J \subset I$ there
exist a function $\varphi = \varphi_J \in C^2(\CR^d)$
diverging to $+\infty$ as $|x|$ tends to $+\infty$, and a positive constant $M_J$
such that
\[
(\A(t)\varphi)(x)\le M_J,\qquad t\in J,\;\,x\in\CR^d.
\]
\end{hyp}

Hypothesis \ref{hyp5} allows to define $G(t,s)$ on a larger class
than $B_{b}(\CR^d)$. Namely, we
show that the right hand side of \eqref{intrep} makes sense for $f=\varphi$,
where $\varphi$ is any of the functions in Hypothesis \ref{hyp5}.

Let us begin with the following fundamental lemma.
If $J\subset I$ is any interval, we set
$$\Lambda _{J}:= \{(t,s)\in J\times J: s\le t\}. $$

\begin{lem}\label{lemma3.5}
Assume that Hypotheses  $\ref{hyp1}(i)(ii)$  and  $\ref{hyp5}$  are satisfied.
Fix a bounded interval $J\subset I$ and let $\varphi=\varphi_J$ be as in
Hypothesis $\ref{hyp5}$.
Then, the function $(t,s,x)\mapsto (G(t,s)\varphi)(x)$ is well defined
and bounded in $\Lambda _{J}\times B_{\varrho} $, for every $\varrho >0$.
\end{lem}
\begin{proof}
We may assume (possibly adding a constant) that $\varphi(x) \geq 0$
for each $x\in \CR^d$.

For every $n\in \CN$ choose a function $\psi_n \in C^{\infty}([0,
+\infty))$ such that
\begin{enumerate}
\item[(i)] $\psi_n(t) = t$ for $t \in [0, n]$,
\item[(ii)] $\psi_n(t) \equiv \mathrm{const.}$ for $ t \geq n+1$,
\item[(iii)] $0 \leq \psi_n' \leq 1$ and $\psi_n'' \leq 0$.
\end{enumerate}
Then, the function $\varphi_n := \psi_n\circ \varphi$ belongs to
$C^2_b(\CR^d)$ and it is constant outside a compact set. By
Lemma \ref{sderivative}, we have
\begin{eqnarray}
\varphi_n(x) &\hs{5} \geq\hs{5} & \varphi_n(x) - (G(t,s)\varphi_n)(x)\nonumber\\
&\hs{5} =\hs{5} & - \int_s^t \int_{\CR^d} (\A(r)\varphi_n)(y) \,
p_{t,r}( x, dy)\,dr\nonumber\\
&\hs{5} =\hs{5} & - \int_s^t \int_{\CR^d} \left \{\psi_n'(\varphi) (\A(r)
\varphi) (y) + \psi_n''(\varphi) \langle Q(r,y) \nabla\varphi
(y) , \nabla\varphi (y) \rangle\right\}
p_{t,r}(x,dy)\,dr \nonumber\\
&\hs{5} \geq\hs{5} & - \int_s^t\int_{\CR^d} \psi_n'(\varphi) (\A(r)
\varphi) (y) \, p_{t,r}(x,dy)\,dr,
\label{chain}
\end{eqnarray}
for any $s,t\in \Lambda$  and any $x\in\CR^d$.
We claim that for each  $s$, $t\in J$, letting  $n\to +\infty$ in \eqref{chain}
we obtain
\[
\varphi (x) \geq - \int_s^t \int_{\CR^d} (\A(r) \varphi) (y) \,
p_{t,r}(x,dy) \, dr=-\int_s^t(G(t,r)\A(r) \varphi) (x)\,dr,
\]
so that, in particular, the above integral is finite.

It is clear that $\lim_{n\to +\infty} \varphi_n(x)= \varphi(x)$ for
each $x\in\CR^d$. Concerning the integral in the right-hand side of \eqref{chain}, we split it
into the sum
\begin{eqnarray}
&\hs{5}& \int_s^t\int_{\CR^d} \psi_n'(\varphi) (\A(r)
\varphi) (y) \, p_{t,r}(x,dy)\,dr\nonumber\\
&\hs{5}=\hs{5}&
-\int_s^t\int_{\CR^d} \psi_n'(\varphi) \left\{M_J-(\A(r)
\varphi) (y)\right\}\, p_{t,r}(x,dy)\,dr\nonumber\\
&\hs{5}&+M_J\int_s^t\int_{\CR^d} \psi_n'(\varphi)\, p_{t,r}(x,dy)\,dr.
\label{may-9}
\end{eqnarray}
Since $\psi_n'(\varphi)(y)$ is increasing in $n$ and converges to $1$
for each $y$,
both integrals in the right-hand side of \eqref{may-9} converge by the monotone convergence theorem.
The claim follows.

Letting  $n\to +\infty$ in  \eqref{chain} yields
\begin{equation}
(G(t,s)\varphi) (x) \leq \varphi (x) + \int_s^t(G(t,r)\A(r)\varphi) (x) \, dr,
\label{sab-marzo}
\end{equation}
and since $(\A(r)\varphi)(y)\leq M_{J}$ for each $y\in \CR^{d}$ and
  $r\in J$,
\begin{equation}
\int_s^tG(t,r)(\A(r)\varphi)(y)\,dr\le M_J(t-s).
\label{883-1}
\end{equation}
Estimates \eqref{sab-marzo} and \eqref{883-1} imply that
\begin{eqnarray*}
(G(t,s)\varphi) (x) \leq \varphi (x) + M_J(t-s),
\end{eqnarray*}
for any $s,t\in J$, with $s\le t$ and any $x\in\CR^d$.
It follows that
\begin{equation}
M_{J,\varrho}:=\sup_{{(t,s,x)\in J\times J\times\CR^d}\atop{s\le t,~|x|\le\varrho}}(G(t,s)\varphi) (x)<+\infty.
\label{Mjrho}
\end{equation}
This completes the proof.
\end{proof}

Having $(G(t,s)\varphi) (x)$ bounded for
$(t,s)\in \Lambda _{J} $, we may prove in the standard
way that for each $r>0$ the family of measures
$\{p_{t,s}(x,dy):~(t,s,x)\in \Lambda _{J}\times\overline{B_r}\}$
is tight. We recall that a family of (probability)
measures $\{\mu_{\alpha}:~\alpha\in {\mathcal F}\}$ is tight, if for any $\varepsilon>0$ there exists
$\varrho>0$ such that $\mu_{\alpha}(\CR^d\setminus B(\varrho))\le\varepsilon$ for any
$\alpha\in {\mathcal F}$.

\begin{lem}
    \label{tightness}
Under the assumptions of Lemma $\ref{lemma3.5}$,
for each bounded interval $J\subset I$ and for each $r>0$ the family of measures
    $\{p_{t,s}(x,dy):~(t,s,x)\in \Lambda _{J}\times\overline{B_r}\}$
    is tight.
\end{lem}
\begin{proof}
Fix $\varepsilon>0$ and consider the function $\varphi=\varphi_J$ in
Hypothesis \ref{hyp5}. As in the proof of Lemma \ref{lemma3.5}, we
assume that $\varphi$ is nonnegative. Since $\varphi$ blows up as $|x|\to +\infty$,
there exists $\varrho >0$ such that
\begin{eqnarray*}
    \varphi(x)\ge \frac{M_{J,r}}{\varepsilon}(\one_{\CR^d\setminus B(\varrho)})(x),\qquad\;\,x\in \CR^d,
    \end{eqnarray*}
where $M_{J,r}$ is given by \eqref{Mjrho}. Then, for $(t,s)\in
\Lambda_{J}$, we have
\begin{eqnarray*}
    p_{t,s}(x,\CR^d\setminus B(\varrho))&\hs{5}=\hs{5}&\int_{\CR^d}\one_{\CR^d\setminus B(\varrho)}(y)p_{t,s}(x,dy)
    \nonumber\\
    &\hs{5}\le\hs{5}& \frac{\varepsilon}{M_{J,r}}\int_{\CR^d}\varphi(y) p_{t,s}(x,dy)\nonumber\\
    &\hs{5}=\hs{5}&
    \frac{\varepsilon}{M_{J,r}}(G(t,s)\varphi)(x)\le\varepsilon,
\end{eqnarray*}
so that
\begin{equation}
    \sup_{(t,s,x)\in \Lambda_{J}\times\overline{B_r}}p_{t,s}(x,\CR^d\setminus B(\varrho))\le
    \varepsilon,
    \label{tight-1}
    \end{equation}
and the statement follows.
\end{proof}

As usual,  tightness yields some convergence result.

\begin{prop}
\label{prop-luca}
Assume that Hypotheses $\ref{hyp1}(i)$-$(ii)$ and $\ref{hyp5}$ are satisfied.
Further, let $\{f_n\}$ be a bounded sequence in $C_b(\CR^d)$, such that $\|f_n\|_{\infty} \leq M$ for each $n\in \CN$ and $f_n$ converges
to $f\in C_b(\CR^d)$ locally uniformly in $\CR^d$.
Then, the
function $G(\cdot,\cdot)f_n$ converges to $G(\cdot,\cdot)f$ locally
uniformly in $\Lambda \times\CR^d$.
\end{prop}
\begin{proof}
Fix any bounded interval $J\subset I$ and any $\varepsilon$, $r>0$. Let   $\varrho$ be such that
\eqref{tight-1} holds, and for $(t,s,x)\in \Lambda_{J}\times \overline{B_{\rho}}$  split $G(t,s)f_n-G(t,s)f$ as
\begin{eqnarray*}
(G(t,s)f_n)(x)-(G(t,s)f)(x)&\hs{5}=\hs{5}&\int_{\CR^d}(f_n(y)-f(y))p_{t,s}(x,dy)\\
&\hs{5}=\hs{5}&\int_{B(\varrho)}(f_n(y)-f(y))p_{t,s}(x,dy)\\
&\hs{5}&+ \int_{\CR^d\setminus B(\varrho)}(f_n(y)-f(y))p_{t,s}(x,dy),
\end{eqnarray*}
so that
\begin{eqnarray*}
|(G(t,s)f_n)(x)-(G(t,s)f)(x)|&\hs{5}\le\hs{5} & \sup_{y\in B(\varrho)}|f_n(y)-f(y)|\int_{\CR^d}p_{t,s}(x,dy)
\\
&\hs{5}&+\left (\sup_{n\in\CN}\|f_n\|_{\infty}+\|f\|_{\infty}\right )\int_{\CR^d\setminus B(\varrho)}p_{t,s}(x,dy)
\\
&\hs{5}\le\hs{5}& \sup_{y\in B(\varrho)}|f_n(y)-f(y)|
\\
&\hs{5}&+ 2M\varepsilon .
\end{eqnarray*}
Fix $n_0\in\CN$ such that
\begin{eqnarray*}
\sup_{y\in B(\varrho)}|f_n(y)-f(y)|\le\varepsilon,\qquad\;\,n\ge n_0.
\end{eqnarray*}
For   $n\geq n_{0}$  we get
\begin{eqnarray*}
&\hs{5}&\sup_{(t,s,x)\in \Lambda_J\times\overline{B_r}}|(G(t,s)f_n)(x)-(G(t,s)f)(x)|
 \\
&\hs{5}=\hs{5}&\varepsilon (1+2M ).
\end{eqnarray*}
Thus, $G(\cdot,\cdot)f_n$ converges to $G(\cdot,\cdot)f$ uniformly in
$\Lambda_{J}\times \overline{B_r}$.
\end{proof}

Now we are ready to prove that $(t,s,x)\mapsto
(G(t,s)f)(x)$ is continuous, for each  $f\in C_b(\CR^d)$.

\begin{thm}
\label{thm-luca}
Under the assumptions of Proposition $\ref{prop-luca}$, the function
$(t,s,x)\mapsto (G(t,s)f)(x)$ is continuous in $\Lambda \times\CR^d $,
for any $f\in C_b(\CR^d)$.
\end{thm}

\begin{proof}
Fix $f\in C_b(\CR^d)$ and let $\{f_n\}\in C^{\infty}_c(\CR^d)$ be a sequence of smooth functions
converging to $f$ locally uniformly in $\CR^d$ and such that
$$
\sup_{n\in\CN}\|f_n\|_{\infty}<+\infty.
$$
By Proposition \ref{prop-luca}, the sequence of functions $(t,s,x)\mapsto
(G(t,s)f _n)(x)$ converges to $ (t,s,x)\mapsto (G(t,s)f)(x)$  locally uniformly.
Therefore, it suffices to show that $ (t,s,x)\mapsto (G(t,s)g)(x)$ is continuous in $\Lambda\times \CR^d$ whenever
$g\in C^{\infty}_c(\CR^d)$. For this purpose, we observe that the classical interior Schauder estimates as in
\cite[Theorem 3.5]{friedman} imply a slightly more general estimate than \eqref{uptobound}, i.e.,
\begin{equation}
\sup_{s\in [a,b]}\|G_n(\cdot,s)g\|_{C^{1+\alpha/2,2+\alpha}([s,s+m]\times B(m))}\le C\|g\|_{C^{2+\alpha}_b(\CR^d)},
\label{dom-0}
\end{equation}
for any $a,b\in I$, $a<b$,  and some positive constant $C$, independent of $n>m$.

Since the sequence of functions $(t,x)\mapsto (G_n( t,s)g)(x)$ converges
to $(t,x)\mapsto (G(t,s)g)(x)$ in $C^{1+\alpha/2,
2+\alpha}([s,s+m]\times B(m))$ for any $s\in [a,b]$, it follows  that
$(t,x)\mapsto (G(t,s)g)(x)\in C^{1+\alpha/2,2+\alpha}([s,s+m]\times B(m))$ for any $s\in [a,b]$ and its
$C^{1+\alpha/2,2+\alpha}$-norm is bounded by
$C\|g\|_{C^{2+\alpha}_b(\CR^d)}$, with the constant $C$ of formula
\eqref{dom-0}.

Fix $(t_0,s_0,x_0)$, $(t,s,x)\in \Lambda \times \CR^d $, with
$s_{0}$, $s\in [a,b]$.
Suppose that $s_0\le s$. Then, $(t,s_0)\in \Lambda $, and
\begin{eqnarray}
|(G(t,s)g)(x)-(G(t_0,s_0)g)(x_0)|&\hs{5}\le\hs{5}&
|(G(t,s)g)(x)-(G(t,s_0)g)(x)|\nonumber\\
&\hs{5}&+|(G(t,s_0)g)(x)-(G(t_0,s_0)g)(x_0)|.
\label{dom-1}
\end{eqnarray}
By \eqref{dom-int} there exists a positive constant $C=C(a,b)$ such that
\begin{eqnarray}
|(G(t,s)g)(x)-(G(t,s_0)g)(x)|\le C|s-s_0|.
\label{dom-2}
\end{eqnarray}
Combining \eqref{dom-1} and \eqref{dom-2} yields
\begin{eqnarray*}
\lim_{{(t,s,x)\to (t_0,s_0,x_0)}\atop{s\ge s_0}}(G(t,s)g)(x)=
(G(t_0,s_0)g)(x_0).
\end{eqnarray*}

Let now assume that $s<s_0$ and split
\begin{eqnarray}
|(G(t,s)g)(x)-(G(t_0,s_0)g)(x_0)|&\hs{5}\le\hs{5}&
|(G(t,s)g)(x)-(G(t_0,s)g)(x_0)|\nonumber\\
&\hs{5}&+|(G(t_0,s)g)(x_0)-(G(t_0,s_0)g)(x_0)|.
\label{dom-3}
\end{eqnarray}
Since $(t,x) \mapsto (G(\cdot,s)g)(x)$ is continuous in
$[s,+\infty)\times\CR^d$, locally
uniformly with respect to $s$,
from \eqref{dom-2} and \eqref{dom-3} we also deduce that
\begin{eqnarray*}
\lim_{{(t,s,x)\to (t_0,s_0,x_0)}\atop{s< s_0}}(G(t,s)g)(x)=
(G(t_0,s_0)g)(x_0).
\end{eqnarray*}
This completes the proof.
\end{proof}


\section{Gradient estimates}
\label{gradient-sect}
In this section we prove both uniform and pointwise gradient estimates.
Besides being interesting in their own right, we will need them in the
next section to prove uniqueness of systems of invariant measures in a
suitable class and convergence results.

Throughout the section we assume that Hypotheses \ref{hyp1} and
\ref{hyp2} hold. Therefore, the bounded classical solution
of problem \eqref{pde}
is such that its first-order spatial derivatives belong to $C^{1+\ahalf , 2 +\alpha}_{\rm loc}((s,+\infty)\times\CR^d)$
(see e.g., \cite[Theorem 3.10]{friedman} and \cite{lieberman}).
We will use this fact in the sequel to
apply a variant of the Bernstein method to get our gradient estimates.

First, we prove uniform gradient estimates.

\begin{thm}\label{uniformgrad}
Let $s\in I$ and
$T > s$. Then, there exist positive constants $C_1,C_2$, depending on
$s$ and $T$, such that:
\begin{enumerate}
\item[(i)]
for every $f \in C^1_b(\CR^d)$ we have
\[ \|\nabla G(t,s)f\|_{\infty} \leq
C_1\|f\|_{C^1_b(\CR^d)}, \qquad\;\,
s<t\leq T;
\]
\item[(ii)]
for every $f \in C_b(\CR^d)$ we have
\[ \|\nabla G(t,s)f\|_{\infty} \leq
\frac{C_2}{\sqrt{t-s}}\|f\|_{\infty}, \qquad\;\, s< t \leq
T.
\]
\end{enumerate}
\end{thm}
\begin{proof}
It suffices to prove the statement for $f\in
C^{2+\alpha}_c(\CR^d)$, since we may approximate an arbitrary $f$ by a
sequence $(f_n) \subset
C^{2+\alpha}_c(\CR^d)$, bounded with respect to the sup-norm and converging to $f$
locally uniformly in $\CR^d$, and Step 3 of
Theorem \ref{thm1} shows that $\nabla G(\cdot,s)f_n$ converge
to $\nabla G(\cdot,s)f$ pointwise in $(s,T] \times \CR^d$.

Let   $k$, $\rho$ be the functions in Hypothesis \ref{hyp2}. Set
\begin{eqnarray*}
k_0:= \sup_{t\in [s,T]}k(t),\qquad\;\, \rho_0:=
\sup_{t\in [s,T]}\rho (t).
\end{eqnarray*}

(i). Let $u_n$ be the unique solution of the Cauchy-Neumann problem \eqref{approx-neumann}, where $n$ is
so large that the support of $f$ is contained in $B_n$. By
Remark \ref{rem2.3}, $u_n$ converges to $u(t,x):=
(G(t,s)f)(x)$  in $C^{1,2}([s,T]\times K)$ as $n\to +\infty$, for any compact set $K\subset \CR^d$.

Define
\[
z_n(t,x) = u_n(t,x)^2 + a|\nabla_x u_n(t,x)|^2,\qquad\;\, (t,x)\in
[s,T]\times B_n.
\]
Then, $z_n$
belongs to $C^{1,2}((s, +\infty )\times B_n)\cap C_b([s,T]\times \overline{B_n})$, for any $s<T$.
Since $B_n$ is convex,
the matrix $D\nu=(D_j\nu_i)$ is positive definite. Moreover, differentiating the
equality
$\frac{\partial u}{\partial \nu} = 0$, one easily verifies that
\begin{eqnarray*}
\sum_{i,j=1}^d \nu_j D_{ij}u D_i u=
-\sum_{i,j=1}^dD_i\nu_j D_i u D_ju\le 0,
\end{eqnarray*}
which, in its turn, implies that the normal derivative of $z_{n}$ on $\partial B_n$ is nonpositive.

We claim that we may choose $a>0$ in such a way that $D_tz_n- \A(t)zn
\leq 0$ for $s<t<T$. Then, the classical maximum principle yields
$|z_n| \leq  \|f\|_{C^1_b(\CR^d)}^2$, i.e.,
\[
u_n(t,x)^2 + a|\nabla u_n(t,x)|^2 \leq \|f\|_{C^1_b(\CR^d)}^2\qquad\;\,
(t,x) \in (s,T) \times B_n.
\]
Letting   $n \to +\infty$, statement (i) follows with $C_1=
a^{-\frac{1}{2}}$.

From now on we  omit the subscript $n$ as well as the dependence
on $t$ and $x$ to simplify notation.
To prove the claim, observe that
\begin{eqnarray}
z_t - A(\cdot)z &\hs{5} =\hs{5} & 2a\langle \nabla_x b\, \nabla_x u,\nabla_x u\rangle - 2
\langle Q \nabla_x u , \nabla_x u \rangle  - 2a \sum_{k=1}^d \langle Q
\nabla_x D_k u, \nabla_x D_k u \rangle\nonumber\\
&\hs{5}& + 2a \sum_{k=1}^d D_ku \cdot \tr \left(D_kQ\cdot
D_x^2u\right).
\label{equation-z}
\end{eqnarray}
Using Hypothesis \ref{hyp2}(iii), we estimate the last term as
follows,
\begin{eqnarray*}
\left|\sum_{k=1}^d D_k u\cdot \tr (D_k Q D_x^2u) \right|
\le \rho_0 \eta
\sum_{k=1}^d|D_ku| \cdot \sum_{i,j=1}^d |D_{ij}u|
\le\rho_0\eta d^{\frac{3}{2}}|\nabla_x u||D^2_xu|.
\end{eqnarray*}
The other terms are easily estimated
 using Hypotheses \ref{hyp1}(ii) and \ref{hyp2}(ii). Eventually, we get
\begin{eqnarray*}
z_t - A(\cdot)z
&\hs{5}\leq\hs{5} & 2(ak_0-\eta ) |\nabla_x u|^2  - 2a\eta \,
|D^2_xu|^2 + 2a\rho_0\eta d^{\frac{3}{2}}\, |\nabla_x u| |D^2_xu|\\[1mm]
&\hs{5} \leq\hs{5} & 2(ak_0-\eta ) |\nabla_x u|^2  - 2a\eta \, |D^2_xu|^2
+ a\eta \left(\rho_0^2d^3 |\nabla_x u|^2+|D^2_xu|^2
\right)\\
&\hs{5} \le\hs{5} & (2ak_0 - 2 \eta + a\eta \rho_0^2d^3) |\nabla_x u|^2.
\end{eqnarray*}
The right hand side is negative, if we choose $a \leq  d^{-3}\rho_0^{-2}$
such
that $2ak_0 \leq \eta_0$.
\medskip
\par
(ii). We proceed similarly to  (i), defining
\[
z_n(t,x) = u_n(t,x)^2 + a(t-s)|\nabla_x u_n|^2,\qquad\;\,(t,x)\in [s,T]\times B_n.
\]
As above, in what follows we omit the subscript $n$ as well as the dependence on $t$ and $x$.

If we proceed as in part (i), we see that $z$ satisfies an equality
similar to  \eqref{equation-z} with $a$ replaced by
$a(t-s)$  and a further addendum $a|\nabla_x u|^2$. Hence,
\[ z_t - A(\cdot)z \leq (2a (T-s)k_0^+ - 2\eta + a(T-s)\eta d^3\rho_0^2+a )|\nabla_x
u |^2,
\]
where $k_0^+=\max\{k_0,0\}$.
The right-hand side is nonpositive, if we choose $a=a_T\leq
(T-s)^{-1}d^{-3}\rho_0^{-2}$ such that $2a(T-s)k_0+a \leq \eta_ 0$.
By the maximum principle  we obtain $z_n\le \|f\|_{\infty}^2$ and
statement (ii) follows, with $C_2=a^{-\frac{1}{2}}$, letting $n\to +\infty$.
\end{proof}

\begin{rem}
{\rm
In the proof of Theorem \ref{uniformgrad} we have chosen to approximate $G(t,s)f$
by solutions of Cauchy-Neumann problems instead of Cauchy-Dirichlet problems as in the first part
of the paper.
Approximation by Cauchy-Dirichlet problems is in fact possible, but it requires stronger conditions on the
coefficients (see e.g., \cite[Section 6.1]{lm} for the autonomous case), that we want to avoid here.
}
\end{rem}

As a  consequence of Theorem \ref{uniformgrad}, our evolution family
enjoys  the strong Feller property.

\begin{cor}\label{strongfeller}
For any $f
\in B_b(\CR^d)$ and any $t>s$ we have $G(t,s)f \in C_b(\CR^d)$.
\end{cor}

\begin{proof}
Let $f \in B_b(\CR^d)$. Then, there exists a bounded sequence $(f_n)
\subset C_b(\CR^d)$ which converges pointwise to $f$ almost everywhere in $\CR^d$.
As a consequence of Theorem \ref{uniformgrad}, for any fixed $s<t$,
the function $t\mapsto G(t,s)f_n $
is Lipschitz continuous with Lipschitz constant independent of $n$.
The statement follows, observing that, by the dominated convergence theorem
and \eqref{intrep},
$G(t,s)f_n$ converges to $G(t,s)f$ pointwise.
\end{proof}

\begin{cor}\label{continuous}
For any
$f\in C^1_b(\CR^d)$ and $s \in I$, the function
$\nabla G(\cdot,s)f$ is continuous in $[s, +\infty) \times
\CR^d$.
\end{cor}

\begin{proof}
We  have to show only continuity at $t=s$. For any $n\in\CN$, let $\varphi \in
C^{\infty}_c(B_n)$ be such that $0 \leq \varphi \leq 1$ and
$\varphi\equiv 1$ in $B_{n-1}$. Put $u(t,x):=(G(t,s)f)(x)$ and $v = \varphi u$. We have
$v_t - A(t)v = \psi$ in $B_n$, where
\[
\psi = - u A(t)\varphi - 2\langle Q \nabla \varphi , \nabla_x u
\rangle.
\]
 From Theorems \ref{thm1} and \ref{uniformgrad}, it follows that the functions
$u$ and $\nabla_x u$ are bounded and continuous in $(s,T]\times\CR^d$, for any $T>s$.
Since $\varphi$ is compactly
supported in $B_n$, $\psi \in C((s,s+1],C_0(B_n))$.
Moreover,  Theorem \ref{uniformgrad} yields that $\|\psi
\|_{\infty}\leq C\|f\|_{C^1_b(\CR^d)}$ for some $C>0$.

Let $\{G_n(t,s)\}$ be the evolution family
associated with problem \eqref{approx-dirichlet}.
By the variation of constants formula (e.g.,
\cite[Proposition 3.2]{acqui}) we have
\begin{equation}
v(t, \cdot) = G_n(t,s)(\varphi f) + \int_s^tG_n(t,\sigma ) \psi (\sigma ) \,
d\sigma, \qquad\;\, s<t < s+1.
\label{representation}
\end{equation}

By classical gradient estimates (\cite[Chapter IV, Theorem 17]{lsu}),
we get
\[ \|\nabla_{x} G_n(t, \sigma )\psi (\sigma )\|_{\infty}
\leq \frac{C_1}{\sqrt{t-\sigma}}\|\psi (\sigma )\|_{\infty} \leq
\frac{C_2}{\sqrt{t-\sigma}} \|f\|_{C^1_b(\CR^d)},
\]
\noindent
for any  $s<\sigma < t \leq s+1$ and some positive constants
$C_1$ and $C_2$.
Hence, we can differentiate \eqref{representation}
obtaining
\[
\nabla_{x} v(t) = \nabla_{x} G_n(t, s)(\varphi f) + \int_s^t \nabla_{x} G_n(t,
\sigma ) \psi (\sigma ) \, d\sigma,\qquad\;\, s<t < s+1.
\]
Therefore, for any $x, x_0 \in B_{n-1}$ we have
\[
| \nabla_x u(t,x) - \nabla f(x_0)| \leq | (\nabla_x G_n(t,s )(\varphi
f))(x) - \nabla f(x_0)|
+ 2C_2\|f\|_{C^1_b(\CR^d)}(t-s)^{\frac{1}{2}},
\]
and this implies that $\nabla G(\cdot,s)f$ is continuous at the point $(s,x_0)$ since the function
$\nabla_x G_n(\cdot,s )(\varphi f)$ is continuous in $\{s\}\times\overline{B_n}$ by classical results.
Since $n$ is arbitrary, the statement follows.
\end{proof}

Next, we prove a pointwise gradient estimate.

\begin{thm}\label{pointwisegrad}
Assume that Hypotheses $\ref{hyp1}$, $\ref{hyp2}(i)(iii)$ and $\ref{hyp3}(i)$ are satisfied.   Then
for every $p\ge p_0$ and any $f\in C^1_b(\CR^d)$  we have
\begin{equation}
|(\nabla G(t,s)f)(x)|^p \leq e^{\sigma_p(t-s)}(G(t,s)|\nabla f|^p)(x), \qquad\;\, t\geq s,\;\,x \in \CR^d,
\label{grad-punt}
\end{equation}
where
\begin{equation}
\sigma_p = p\, \sup_{(t,x)\in I\times\CR^d}\left(r(t,x) + \frac{d^3(\rho(t))^2\eta(t,x)}{4\min\{p_0-1,1\}}\right).
\label{sigma-p-1}
\end{equation}

Similarly, under Hypotheses $\ref{hyp1}$, $\ref{hyp2}(i)$ and $\ref{hyp3}(ii)$,  estimate \eqref{grad-punt}
holds true for any $p\in (1,+\infty)$, with
\begin{equation}
\sigma_p=p\left (k_0+\frac{d^3\rho_{0}^{2}}{4\min\{p-1,1\}}\right ).
\label{sigma-p-2}
\end{equation}

Moreover, if the coefficients $q_{ij}$ $(i,j=1,\ldots,d)$ do not depend on $x$ and $r\le\eta$
in $I\times\CR^d$,
then \eqref{sigma-p-2} holds true
for $p=1$ too, provided Hypothesis $\ref{hyp3}(ii)$ is satisfied. In such a case, $\sigma_1=k_0$.
\end{thm}

\begin{proof}
To prove the first part of the statement, fix $s \in I$ and $\eps > 0$. Set $u(t,x): = (G(t,s)f)(x)$ and define
\[
w(t,x) = \left(|\nabla_x u(t,x)|^2 + \eps\right)^{\frac{p}{2}}, \qquad\;\, t \geq s,\;\,
x \in \CR^d.
\]
By Corollary \ref{continuous},   $w \in C_b([s, T] \times
\CR^d)$ for all $T>s$, and moreover by \cite[Theorem
3.10]{friedman},  $w \in C^{1,2}((s,T) \times \CR^d )$.
A straightforward computation shows that
\[
w_t - \A(t)w = f_1 + f_2 + f_3,
\]
where
\begin{eqnarray*}
f_1 &\hs{5} =\hs{5} & p( |\nabla_x u|^2 + \eps )^{\frac{p}{2}-1} \left(
\langle \nabla_x b\, \nabla_x u , \nabla_x u \rangle - \sum_{k=1}^d \langle
Q\nabla_x D_k u , \nabla_x D_k u\rangle\right),
\\
f_2 &\hs{5} =\hs{5} & p( |\nabla_x u|^2 + \eps )^{\frac{p}{2}-1}\sum_{k=1}^d
D_ku \cdot\tr (D_k Q \cdot D^2_xu),\\
f_3 &\hs{5} =\hs{5} & -p(p-2)(|\nabla_x u|^2 + \eps )^{\frac{p}{2}-2} \langle Q
D^2_xu\, \nabla_x u, D^2_xu\, \nabla_x u\rangle.
\end{eqnarray*}
Using Hypotheses \ref{hyp1}(ii), \ref{hyp2}(i)(iii) and \ref{hyp3}(i), we   estimate $f_1$  as in the proof of
Theorem \ref{uniformgrad}, getting
\begin{eqnarray}
f_1
&\hs{5}\le\hs{5}& p(|\nabla_x u|^2 +\eps)^{\frac{p}{2}-1} \cdot \left (
r\, |\nabla_x u|^2 -\sum_{k=1}^d\langle Q\nabla_x D_ku,\nabla_x D_ku\rangle \right )
\label{first}\\
\nonumber
&\hs{5}\le\hs{5}&  p(|\nabla_x u|^2 +
\eps)^{\frac{p}{2}-1} \cdot ( r\, |\nabla_x u|^2 - \eta\, |D_x^2u|^2 ).
\label{second}
\end{eqnarray}
Moreover, for every $c>0$ we have
\begin{equation}
f_2
\le
p(|\nabla_x u|^2 + \eps)^{\frac{p}{2}-1}\eta\left ( c \, |D_x^2u|^2 + \frac{d^3\rho^2}{4c}|\nabla_x u|^2\right ).
\label{third}
\end{equation}

Concerning $f_3$, we have
\begin{eqnarray}
\langle Q D_x^2u \nabla_x u , D_x^2 u\nabla_x u\rangle
&\hs{6}  = \hs{6} &  |Q^{1/2} D_x^2u \nabla_x u|^2
\nonumber\\
&\hs{6} \leq\hs{6} & \|Q^{1/2} D_x^2u\|^2 \, |\nabla_x u|^2
\nonumber\\
&\hs{6}= \hs{6} &
|\nabla_x u|^2 \sum_{k=1}^d \langle Q \nabla_x D_k u , \nabla_x
D_k u \rangle.
\label{star-2}
\end{eqnarray}

Now we distinguish between two cases.
\medskip
\par
{\em Case $1$: $p \ge\max\{p_0,2\}$}.
Since $p(p-2) \geq 0$,  the uniform ellipticity assumption
implies   $f_3 \leq 0$. Using \eqref{first} and \eqref{third} with
$c=1$,
 we obtain
\begin{eqnarray}
w_t - A(\cdot)w
&\hs{5}\le\hs{5}& \sigma_p (|\nabla_x u|^2 + \eps )^{\frac{p}{2}-1}|\nabla_x u|^2\nonumber\\
&\hs{5}=\hs{5}&\sigma_p(|\nabla_x u|^2 + \eps )^{\frac{p}{2}}-\sigma_p\varepsilon
(|\nabla_x u|^2 + \eps )^{\frac{p}{2}-1}.
\label{montecarlo}
\end{eqnarray}
Now, observing that
\begin{eqnarray*}
a^{\frac{p}{2}-1}\le \left (1-\frac{2}{p}\right )a^{\frac{p}{2}}+\frac{p}{2},\qquad\;\,a>0,
\end{eqnarray*}
from \eqref{montecarlo} we deduce that
\begin{eqnarray*}
w_t-A(\cdot)w\le \sigma_{p,\varepsilon}(w- \delta_{\eps}),
\end{eqnarray*}
where,
\begin{eqnarray*}
\sigma_{p,\varepsilon}=
\left\{
\begin{array}{ll}
\sigma_p, &\mbox{ if }\, \sigma_p\ge 0,\\[2mm]
\sigma_p\left \{1-\left (1-\frac{2}{p}\right )\varepsilon\right\}, &\mbox{ if }\, \sigma_p< 0,
\end{array}
\right.
\qquad\;\,
 \delta_{\eps} = \left\{\begin{array}{cl}
0, & \mbox{if}\,\, \sigma_p \geq 0,\\[1mm]
\frac{p}{2}\eps, & \mbox{if}\,\, \sigma_{p} < 0.
\end{array}\right.
\end{eqnarray*}
and $\sigma_{p,\varepsilon}$  is given by \eqref{sigma-p-1}.

{\em Case $2$: $p_0<2$ and $1<p<2$.}
In this case, $-p(p-2)$ is positive. Hence, we may use
\eqref{star-2} to estimate $f_3$.
Together with   estimates \eqref{second} and \eqref{star-2} (with $c=p-1$),   we obtain
\begin{eqnarray*}
w_t - A(\cdot)w&\hs{5} \le\hs{5} & p(|\nabla_x u|^2 + \eps)^{\frac{p}{2}-1}\\
&\hs{5}&\qquad\times
\bigg ( (2-p)\sum_{k=1}^d\langle Q\nabla_x D_ku,\nabla_x D_ku\rangle
+ \frac{d^3\rho^2}{4(p-1)}\eta|\nabla_x u|^2 \\
&\hs{5}&\qquad\qquad\;\, + r
|\nabla_x u|^2 - \eta\sum_{k=1}^d\langle Q\nabla_x D_ku,\nabla_x D_ku\rangle+(p-1)\eta|D^2_xu|^2\bigg )\\
&\hs{5}\le\hs{5}& \sigma_{p,\varepsilon} (w- \delta_{\eps}).
\end{eqnarray*}
Here,
\begin{eqnarray*}
 \delta_{\eps} = \left\{\begin{array}{cl}
0, & \mbox{if}\,\, \sigma_p \geq 0,\\[1mm]
\eps^{\frac{p}{2}}, & \mbox{if}\,\, \sigma_{p} < 0,
\end{array}\right.
\end{eqnarray*}
 and $\sigma_{p,\varepsilon}:=\sigma_p$ is given by \eqref{sigma-p-1}.

Now the procedure is the same in the two cases. Setting  $v = w
- \delta_{\eps}$   we have $v_t - A(\cdot)v \leq \sigma_{p,\varepsilon} v$. On the
other hand, the function
\[ z(t) = e^{\sigma_{p,\varepsilon}(t-s)}G(t,s) (|\nabla f|^2 + \eps)^{\frac{p}{2}},\qquad\;\,t>s,
\]
satisfies $z_t - \A(t)z = \sigma_{p,\varepsilon}z$. Thus,
\[ \left\{ \begin{array}{rcll}
 (v-z)_t - (\A(t) + \sigma_{p,\varepsilon})(v-z) &\hs{5} \leq\hs{5} & 0,\quad\; &t\in (s,+\infty),\\[1.5mm]
(v-z)(s) &\hs{5}=\hs{5}& - \delta_{\eps}.
\end{array}
\right.
\]
Theorem \ref{maximum} implies   $v \leq z$. Letting $\eps
\to 0^+$, the statement follows by Proposition \ref{properties}.

In the case that Hypothesis \eqref{hyp3}(i) is replaced by \eqref{hyp3}(ii), the functions    $f_1$,  $f_2$ are estimated as follows:
\begin{eqnarray}
f_1
&\hs{5}\le\hs{5}& p(|\nabla_x u|^2 +\eps)^{\frac{p}{2}-1} \cdot \left (
k_0\, |\nabla_x u|^2 -\sum_{k=1}^d\langle Q\nabla_x D_ku,\nabla_x D_ku\rangle \right )
\label{first-1}\\
&\hs{5}\le\hs{5}&  p(|\nabla_x u|^2 +
\eps)^{\frac{p}{2}-1} \cdot (k_0\, |\nabla_x u|^2 - \eta\, |D_x^2u|^2 ),
\label{second-1}
\end{eqnarray}
\begin{equation}
f_2
\le
p(|\nabla_x u|^2 + \eps)^{\frac{p}{2}-1}\left ( c\eta \, |D_x^2u|^2 + \frac{d^3\rho^2_0}{4c}|\nabla_x u|^2\right ),
\label{third-1}
\end{equation}
for any  $c>0$.
Then, estimate \eqref{grad-punt} with $p\in (1,+\infty)$ (and with $p=1$ too, if the diffusion coefficients
are constant with respect to $x$), follows arguing as above.
\end{proof}

\begin{cor}\label{uniforminfinity}
Under the hypotheses of Theorem $\ref{pointwisegrad}$, there exists a constant $C$ such that
\[
\|\nabla G(t,s)f\|_{\infty} \leq C \cdot e^{\frac{\sigma_p}{p}(t-s)}\|f\|_{\infty},
\quad f \in C_b(\CR^d),\;  s\in I,  \;t \geq s+1,
\]
for every $p\ge p_0$  if Hypothesis $\ref{hyp3}(i)$ is
satisfied,  and for every $p>1$ if Hypothesis $\ref{hyp3}(ii)$ is satisfied.
\end{cor}

\begin{proof}
By Theorem \ref{uniformgrad}, for any $f \in C_b(\CR^d)$ the function $G(s+1,s)f$ is in $C^1_b(\CR^d)$,
and its $C^1$-norm does not exceed $C_1\|f\|_{\infty}$ for
some $C_1>0$, independent of $f$.
If $t>s+1$, we have, by Theorem \ref{pointwisegrad},
\begin{eqnarray*}
| (\nabla G(t,s)f) (x)|^p &\hs{5}=\hs{5}& |(\nabla G(t,s+1)G(s+1,s)f)(x)|\\[1mm]
&\hs{5} \leq\hs{5}& e^{\sigma_p(t-(s+1))} (G(t,s+1) |\nabla G(s+1,s)f|^p)(x)\\[1mm]
&\hs{5} \leq\hs{5} & e^{\sigma_p(t-(s+1))} \|\nabla G(s+1,s)f\|_{\infty}^p.
\end{eqnarray*}
Thus,
\[ \|\nabla G(t,s)f\|_{\infty}^p \leq e^{\sigma_p(t-(s+1))} \|\nabla
G(s+1,s)f\|_{\infty}^p,
\]
and the statement follows.
\end{proof}


\section{Evolution systems of measures}
\label{evol}
\begin{defn}
\label{def-syst}
Let $\{U(t,s)\}$ be an evolution family of bounded operators on $B_b(\CR^d)$.
A family $(\nu_t)$ of probability measures on $\CR^d$ is an {\it
evolution system of measures} for $\{U(t,s)\}$ if, for every $f \in
B_b(\CR^d)$ and every $s < t$, we have
\begin{equation}
\int_{\CR^d} U(t,s)f \, d\nu_t = \int_{\CR^d} f \, d\nu_s.
\label{invar}
\end{equation}
\end{defn}
Formula \eqref{invar} may be rewritten as  $U^*(t,s)\nu_t = \nu_s$.
It  implies that, if we know a single measure $\nu_{t_0}$ of an evolution
system of measures for $\{U(t,s)\}$, then we know all the measures $\nu_t$ for $t\leq t_0$. In particular,
an evolution system of measures is uniquely determined by its tail $(\nu_t)_{t\geq t_0}$.
\medskip
\par
In this section we give sufficient conditions for the existence of an evolution system
($\mu_t$) of measures associated with the evolution family $\{G(t,s)\}$ and we study the main properties of
($\mu_t$).
As a first step, we note that, for our evolution family $\{G(t,s)\}$, evolution systems of measures
necessarily consist of measures which are equivalent to the Lebesgue measure.

\begin{prop}\label{regular}
If $(\mu_t)$ is an evolution system of measures for $\{G(t,s)\}$
then $(\mu_t)$ is equivalent to the Lebesgue measure.
\end{prop}

\begin{proof}
For each $A \in B(\CR^d)$ and $t\in I$ we have
\[
\mu_t(A) = \int_{\CR^d}(G(t+1,t)\one_A)(x) \, \mu_{t+1}(dx).
\]
By Corollaries \ref{irreducibility} and \ref{strongfeller}, if the Lebesgue measure $|A|$ of $A$ is positive then
$(G(t+1,t)\one_A)(x)$ is positive for each $x\in \CR^d$; therefore $\mu_t(A)>0$. On the other hand, by Proposition
\ref{transition}(ii), if $|A|=0$ then $G(t+1,t)\one_A = 0$, hence $\mu_t(A) =0$.
\end{proof}

To prove  existence of evolution systems of measures we use  a procedure similar to the
Krylov-Bogoliubov Theorem which states that, in the autonomous case, existence of an invariant measure is equivalent to the tightness of a certain set of probability measures.
In our case, the corresponding tightness property is proved under Hypothesis  $\ref{hyp4}$, through the Prokhorov Theorem. It states that a set $\{P_{\alpha}: \alpha\in{\mathcal F}\}$
of probability measures is tight if and only if,
for any sequence $(\alpha_n)$ in $\mathcal{F}$,
there exists a subsequence $\alpha_{n_k}$ such that
$P_{\alpha_{n_k}}$ converges to some probability measure $P$ in the following sense:
\[
\lim_{k\to +\infty} \int_{\CR^d} f(y) \, P_{\alpha_k}(dy) = \int_{\CR^d}
f(y) \, P(dy), \qquad\;\, f\in C_b(\CR^d).
\]

\begin{lem}\label{khashminskii}
Assume that Hypotheses $\ref{hyp1}$ and $\ref{hyp4}$ are satisfied.
Then, $G(t,s)\varphi$ is well defined for any $t_0\le s\le t\in I$. Moreover,
for any fixed $x\in\CR^d$, the function $(t,s)\mapsto (G(t,s)\varphi)(x)$ is bounded
in $\Lambda = \{(t,s)\in I\times I: t_0\le s\le t\}$.
\end{lem}

\begin{proof}
Lemma \ref{lemma3.5} implies that $G(t,s)\varphi $
is well defined for   $(t,s)\in \Lambda $ with $t_0\le s $ and the function $(t,s,x)\mapsto
(G(t,s)\varphi)(x)$ is locally bounded.
To complete the proof, we fix $t>t_0$ and $x\in\CR^d$, and consider the function
$g$  defined in $[t_0,t]$ by $g(s):=(G(t,s)\varphi)(x)$. $g$ is measurable, because $(G(t,s)\varphi)(x)$
is the pointwise limit of the functions $(G(t,s)\varphi_n)(x)$ in the proof of Lemma \ref{lemma3.5},
that are continuous with respect to $s$.
The procedure of Lemma \ref{lemma3.5} yields
\begin{equation}
g(r)- g(s) \geq \int_s^{r}(cg(\sigma)-a)\, d\sigma,\qquad\;\, t_0\le s\le r\le t.
\label{g}
\end{equation}
We claim that \eqref{g} implies
\begin{equation}
 g(s) \leq \bigg( g(t)-\frac{a}{c}\bigg) e^{c(s-t)} + \frac{a}{c},\qquad\;\, t_0\le s\le  t.
\label{g0}
\end{equation}
Indeed, for any fixed $s\geq t_0 $, the function $\Phi$ defined by
$$\Phi(r): = \bigg( g(s) -\frac{a}{c} + \int_s^{r}(cg(\sigma)-a)\, d\sigma\bigg) e^{-cr}, \qquad\;\, r\leq s\leq t$$
is continuous in $[s,t]$ and therein weakly differentiable with $\Phi'(r)\geq 0$  a.e., so that it is nondecreasing,
and $\Phi(s)\leq \Phi(t)$ implies \eqref{g0}.
From \eqref{g0} we obtain $(G(t,s)\varphi)(x) \leq \varphi(x) + a/c$, and the statement follows.
\end{proof}

\begin{thm}\label{krylov}
Assume that Hypotheses $\ref{hyp1}$ and $\ref{hyp4}$ are satisfied.
Then, there exist an evolution system $(\mu_t)$ of measures for $\{G(t,s)\}$
and a constant $M \geq 0$ such that
\begin{equation}
\int_{\CR^d} \varphi (y) \, \mu_t(dy) \leq M,\qquad\;\,t\ge t_0.
\label{merc-sera}
\end{equation}
\end{thm}
\begin{proof}
Fix $s\in I$ and $x_0\in\CR^d$.
For any $t> s$,   define  the measure $\mu_{t,s}$ by
\[
\mu_{t,s}(A) := \frac{1}{t-s}\int_s^t p_{\tau , s}(x_0,A) \, d\tau
= \frac{1}{t-s}\int_s^t (G(\tau , s)\one_A)(x_0)\, d\tau.
\]
Lemma  \ref{khashminskii} implies that  the family $(\mu_{t,s})_{t> s\geq t_0}$
is tight, through the same proof of Lemma \ref{tightness}.
The  Prokhorov Theorem and a diagonal argument yield   existence of a
sequence $t_k$ diverging to $+ \infty$ and of probability measures
$\mu_n$ ($n\in\CN$, $n>t_0$) such that
$\mu_{t_k , n} \weak^* \mu_n$. To define $\mu_s$ also for noninteger $s$, we show preliminarly that
 $G^*(n,m)\mu_n = \mu_m$ for $m < n$. Indeed, for each  $A \in \mathcal{B}(\CR^d)$ we have
\begin{eqnarray*}
G^*(n,m)\mu_n(A) &\hs{5}=\hs{5}& \int_{\CR^d}\one_A(y) \,
G^*(n,m)\mu_n(dy)\\
&\hs{5} =\hs{5}& \int_{\CR^d}(G(n,m)\one_A)(y) \, \mu_n(dy)\\
&\hs{5} =\hs{5} & \lim_{k\to +\infty} \frac{1}{t_k-n}\int_n^{t_k}
(G(\tau , n)G(n,m)\one_A)(x) \, d\tau\\
&\hs{5} =\hs{5} & \lim_{k\to +\infty} \frac{1}{t_k-n}\int_n^{t_k}(G(\tau ,
m)\one_A)(x) \, d\tau\\
&\hs{5} =\hs{5} & \lim_{k \to +\infty} \frac{1}{t_k-m}\int_m^{t_k} (G(\tau , m)\one_A)(x) \, d\tau\\
&\hs{5} =\hs{5} & \mu_m(A).
\end{eqnarray*}
Thus, we can extend the definition of the measures $\mu_s$ to any $s\in I$, by
setting
 $\mu_s:=G^*(n,s)\mu_n$ where $n$ is any positive integer greater than $s$.
Since $G^*(n,s) = G^*(m,s)G^*(n,m)$,   this definition is independent of $n$.
It is immediate to check that $(\mu_t)$ is an evolution system of measures for $\{G(t,s)\}$.

To complete the proof, we observe that, since $(G(t,s)\varphi)(x_0)$ is bounded in $t\ge s\ge t_0$, then each integral
\[
\frac{1}{t-s}\int_s^t(G(\tau,s)\varphi)(x_0)d\tau=\int_{\CR^d}\varphi(y)\mu_{t,s}(dy)
\]
is bounded for $t>s\ge t_0$ by the same constant. Letting $t\to +\infty$, we get \eqref{merc-sera}.
\end{proof}

\begin{rem}
It should be noted that the evolution system of measures constructed
in Theorem \ref{krylov} could still depend on $x_0$. Indeed, in general,
evolution system of measures are not unique.
In \cite[Lemma 2.2]{gl1} it is proved that the
evolution family associated with the operators
\[
(\A(t)u)(x) = \half \Delta u( x) + \langle B(t) x , \nabla u( x) \rangle,
\]
admits infinitely many evolution systems of measures. However,
uniqueness may be achieved among all systems of measures which have finite
moments of order $p$ for some $p> 0$ with a certain asymptotic behaviour.
\end{rem}

In the following, if $(\mu_t)$ is a family of probability measures
on $\CR^d$, we denote by
\[
\mu_t (p) := \int_{\CR^d} |x|^p \, \mu_t(dx ),
\]
the $p$-th moment function.
We note that, if $\varphi (x) = |x|^p$ satisfies Hypothesis \ref{hyp4}, then Theorem
\ref{krylov} implies that $\{G(t,s)\}$ admits an evolution system of measures
$(\mu_t)$ such that $\mu_t (p) = O (1)$ as $t\to +\infty$,
i.e. there exists $t_0\in I$ such that the $p$-th
moments of $\mu_t$ exist and are uniformly bounded
for any $t\geq t_0$.

Let us see the connection between  evolution systems of measures and
asymptotic behaviour of solutions to problem \eqref{pde}.
We
assume that there exists a {\it negative} constant $\omega$ such
that, for large $t-s$, we have
\[
\|\nabla G(t,s)f\|_{\infty} \leq e^{\omega (t-s)}\|f\|_{\infty},
\qquad\;\, f \in C_b(\CR^d).
\]
A sufficient condition for this may be obtained from Corollary
\ref{uniforminfinity}.

\begin{thm}\label{asymptotics}
Assume that there exists $\omega < 0$ such that
\begin{equation}
\|\nabla G(t,s)f\|_{\infty} \leq Ce^{\omega (t-s)}\|f\|_{\infty},
\label{laura}
\end{equation}
for all $t \geq s+1$, all $f \in C_b(\CR^d)$ and some positive constant $C$. Further, assume that
$\{G(t,s)\}$ admits an evolution system of measures $(\mu_t)$ such that, for some $p>0$,
\[
\lim_{t\to +\infty} \mu_t (p)  e^{ \omega pt} =0.
\]
Then,
\[
\lim_{t\to +\infty}(G(t,s)f)(x)=\int_{\CR^d} f(y) \, \mu_s(dy),\qquad\;\,x\in\CR^d,
\]
for all $s\in I$ and $f \in C_b(\CR^d)$. If $I = \CR$,
then, we also have
\[
\lim_{s\to -\infty}\left ((G(t,s)f)(x) - \int_{\CR^d} f(y) \, \mu_s(dy)\right )=0,\qquad\;\,x\in\CR^d.
\]
In both cases the convergence is uniform in the compact sets of $\CR^d$.
\end{thm}

\begin{proof}
Without loss of generality, we may assume that $p<1$. We have
\begin{eqnarray*}
(G(t,s)f)(x) - \int_{\CR^d} f(y) \, \mu_s(dy) &\hs{5}=\hs{5}&  (G(t,s)f)(x) -
\int_{\CR^d}(G(t,s)f)(y)\, \mu_t(dy)\\
&\hs{5}=\hs{5}& \int_{\CR^d} \{(G(t,s) f)(x) - (G(t,s)f)(y)\} \, \mu_t(dy).
\end{eqnarray*}
Splitting
\begin{eqnarray*}
&\hs{5}&|(G(t,s)f)(x) - (G(t,s)f)(y)|\nonumber\\
&\hs{5}=\hs{5}&|(G(t,s)f)(x) - (G(t,s)f)(y)|^{1-p}
|(G(t,s)f)(x) - (G(t,s)f)(y)|^p
\end{eqnarray*}
 and using the mean value theorem and
\eqref{laura}, we get
\begin{eqnarray*}
|(G(t,s)f)(x) - (G(t,s)f)(y)|
\le 2C^p\|f\|_{\infty} e^{p\omega (t-s)} |x-y|^p.
\end{eqnarray*}
Hence, we have:
\begin{eqnarray}
\left |(G(t,s)f)(x) - \int_{\CR^d} f(y) \, \mu_s(dy)\right | &\hs{5} \leq\hs{5} & 2C^p\|f\|_{\infty} e^{p\omega
(t-s)} \int_{\CR^d} |x-y|^p \, \mu_t(dy)\nonumber\\
&\hs{5}\leq\hs{5}& 2C^p\|f\|_{\infty}e^{p\omega (t-s)} \left( |x|^p +  \int_{\CR^d}|y|^p\,
\mu_t(dy)\right),
\label{asympt-behav}
\end{eqnarray}
and the right-hand side   vanishes as $t\to +\infty$ (and also as $s\to -\infty$, if $I=\CR$),
uniformly for $x$ in compact sets.
\end{proof}

\begin{cor}\label{uniquemeas}
Under the hypothesis of Theorem $\ref{asymptotics}$, there exists at most one evolution system of measures
$(\mu_t)$ such that $\lim_{t\to+\infty} \mu_t (p)  e^{ \omega pt} =0 $ for some $p>0$.
\end{cor}
\begin{proof}
Let $(\mu_t)$, $(\nu_t)$ be two evolution system of measures with the above property. By Theorem \ref{asymptotics}, for each $f\in C_b(\CR^d)$ and $s\in I$ we have
$$\int_{\CR^d} f(y) \, \mu_s(dy) =\int_{\CR^d} f(y) \, \nu_s(dy),$$
since both integrals coincide with $ \lim_{t\to +\infty}(G(t,s)f)(0) $. The statement follows.
\end{proof}


\section{Evolution semigroups in $L^p$ spaces with respect to
invariant measures}
\label{sect-5}

In this section we assume that $I = \CR$, and that  Hypotheses \ref{hyp1} and \ref{hyp4} are satisfied.

Let us define the evolution semigroup $\{\T(t)\}$
associated  with the
evolution family $\{G(t,s)\}$ on the space $C_b(\CR^{d+1})$ by
\[
(\T (t)f)(s,x) = (G(s,s-t)f(s-t,\cdot))(x),\qquad\;\,(s,x)\in\CR^{d+1},\;\,t\ge 0.
\]

\begin{prop}
\label{prop6.1}
The family of operators $\{\T(t): \,t\geq 0\}$ is a semigroup of positive contractions in $C_b(\CR^{d+1})$.
Moreover, $\T(t)f$ tends to $f$ locally uniformly in $\CR^{d+1}$ as $t\to 0^+$, for any $f\in C_b(\CR^{d+1})$.
\end{prop}
\begin{proof}
As a first step we prove that, for any $t>0$, the operator $\T(t)$ maps $C_b(\CR^{d+1})$ into itself.
By Theorem \ref{thm1}, we know that
\begin{eqnarray*}
\sup_{(s,x)\in\CR^{d+1}}|(\T(t)f)(s,x)|\le \|f\|_{C_b(\CR^{d+1})},\qquad\;\,t\ge 0.
\end{eqnarray*}
Let us now fix $(s_0,x_0)$ in $\CR^{d+1}$ and observe that
\begin{eqnarray}
&\hs{5}&|(\T(t)f)(s,x)-(\T(t)f)(s_0,x_0)|\nonumber\\
&\hs{5}=\hs{5}&
|(G(s,s-t)f(s,\cdot))(x)-(G(s_0,s_0-t)f(s_0,\cdot))(x_0)|\nonumber\\
&\hs{5}\le\hs{5}& |(G(s,s-t)f(s,\cdot))(x)-(G(s,s-t)f(s_0,\cdot))(x)|\nonumber\\
&\hs{5}&+|(G(s,s-t)f(s_0,\cdot))(x)-(G(s_0,s_0-t)f(s_0,\cdot))(x_0)|.
\label{cont-semigruppone}
\end{eqnarray}
By Proposition \ref{prop-luca},
$$\lim_{s\to s_0}\sup_{(r,x)\in [s_0-\delta,s_0+\delta]\times \{x_0+B(\delta)\}}
|(G(r,r-t)f(s,\cdot))(x)-(G(r,r-t)f(s_0,\cdot))(x)|=0,
$$
for any $\delta>0$.
Therefore,
the first term in the right-hand side of \eqref{cont-semigruppone} converges to $0$ as $(s,x)$ tends to
$(s_0,x_0)$. Similarly, by Theorem \ref{thm-luca}, the function $(p,r,x)\mapsto (G(p,r)f)(x)$ is continuous
in $\{(p,r,x)\in\CR^{d+2}: r\le p\}$. Hence, also the second term   tends to $0$ as $(s,x)$ tends to $(s_0,x_0)$.
This shows that $\T(t)f\in C_b(\CR^{d+1})$.

The semigroup property follows easily since $\{G(t,s)\}$ is an evolution family.
Indeed, for any $t_1<t_2$, it holds that
\begin{eqnarray*}
&\hs{5}&(\T(t_2)\T(t_1)f)(s,x)\\
&\hs{5}=\hs{5}&(G(s,s-t_2)\T(t_1)f(s-t_2,\cdot))(x)\\
&\hs{5}=\hs{5}&\left (G(s,s-t_2)G(s-t_2,s-t_2-t_1)f(s-t_1-t_2,\cdot)\right )(x)\\
&\hs{5}=\hs{5}&\left (G(s,s-t_2-t_1)f(s-t_1-t_2,\cdot)\right )(x)\\
&\hs{5}=\hs{5}&(\T(t_1+t_2)f)(s,x),
\end{eqnarray*}
for any $(s,x)\in\CR^{d+1}$.

The positivity of $\T(t)$ follows from the positivity of the evolution
family $\{G(t,s)\}$.

Finally, the fact  that $\T(t)f$ converges to $f$ locally uniformly in $\CR^{d+1}$
as $t\to 0^+$, is an immediate consequence of the
continuity of the function $(p,r,x)\mapsto (G(p,r)f)(x)$ in $\{(p,r,x)\in\CR^{d+2}:p\ge r\}$
and Proposition \ref{prop-luca}.
\end{proof}

\begin{rem}
{\rm Since \mbox{$G(s,s-t) \one = \one$} for each $t>s$, if $f=f(s)$ depends only on time then
$(\T(t)f)(s,x) = f(s-t)$, i.e.
$\T(t)$ acts as a translation semigroup. Therefore, $\T(t)$ cannot have any smoothing or
summability improving property in the $s$ variable.
In particular, it is not strong Feller and not hypercontractive.
}
\end{rem}

Now, let $(\mu_t)$ be an evolution system of measures for $\{G(t,s)\}$.
Note that the function $s\mapsto \mu_s(A)$ is measurable in $I$ for any Borel set $A$.
Indeed, by Lemma \ref{sderivative}, the function $s\mapsto (G(t,s)f)(x)$
is bounded and continuous in $(-\infty,t)$, for any $x\in\CR^d$ and any $f\in C_0(\CR^d)$.
Hence, the function
\begin{eqnarray*}
s\mapsto\int_{\CR^d}(G(t,s)f)(x)\,\mu_t(dx),
\end{eqnarray*}
is continuous as well in $(-\infty,t)$.
Since
\begin{eqnarray*}
\mu_s(A)=\int_{\CR^d}(G(t,s)\one_A)(x)\,\mu_t(dx),
\end{eqnarray*}
and \mbox{$\one_A$} is the pointwise limit of a sequence  $(f_n)\subset C_0(\CR^d)$, bounded
with respect to the sup-norm, by dominated convergence, the measurability of the function
$s\mapsto\mu_s(A)$ follows.
Therefore, we can define
\[
\nu ( J \times K) := \int_J \mu_t(K) \, dt,
\]
for Borel sets $J \subset \CR$ and  $K \subset
\CR^d$. Of course, $\nu$ may be uniquely extended in a standard way to a
measure on $\mathcal{B}(\CR^{d+1})$.\\

In the following, we denote by $\mathcal{G}$ the differential
operator
\begin{equation}
\mathcal{G}u(t,x) = \A(t)u(t,x) -  u_t(t,x),\qquad\;\,(t,x)\in\CR^{d+1}.
\label{operat-G}
\end{equation}

We state a preliminary lemma about $\T(t)$ and $\mathcal{G}$.

\begin{lem}\label{behoncont}
\begin{enumerate}
\item[(i)]
For all $\varphi \in C_c(\CR, C_b(\CR^d))$ and all $t \geq 0$
we have
\[
\int_{\CR^{d+1}} \T(t) \varphi \, d\nu = \int_{\CR^{d+1}} \varphi \, d\nu.
\]
\item[(ii)]
For $\varphi \in C^{1,2}_c(\CR^{d+1})$ we have
\begin{equation}
\int_{\CR^{d+1}} \mathcal{G}\varphi \, d\nu = 0.
\label{inf-inv}
\end{equation}
\end{enumerate}
\end{lem}

\begin{proof}
(i). We have
\begin{eqnarray*}
\int_{\CR^{d+1}}\T (t)\varphi\, d\nu &\hs{5} =\hs{5} & \int_{\CR}\int_{\CR^d} \left(
G(s,s-t)\varphi (s-t , \cdot)\right) (x) \, \mu_{s}(dx)\,ds\\
&\hs{5} =\hs{5} & \int_{\CR} \int_{\CR^d} \varphi (s-t , x) \, \mu_{s-t}(dx)\,ds\\
&\hs{5} =\hs{5} & \int_{\CR}\int_{\CR^d} \varphi (r , x) \, \mu_r(dx)\,dr\\
&\hs{5} =\hs{5} & \int_{\CR^{d+1}} \varphi \, d\nu.
\end{eqnarray*}
\medskip
\par
(ii). By part (i) we obtain
\[
\lim_{h \to 0^+} \int_{\CR^{d+1}} \frac{\T (h)\varphi - \varphi
}{h} \, d\nu = 0.
\]
Now we show that
\[
\lim_{h \to 0^+} \int_{\CR^{d+1}} \frac{\T (h)\varphi - \varphi
}{h} \, d\nu =  \int_{\CR^{d+1}} \mathcal{G}\varphi \, d\nu .
\]
For this purpose, let $a,b\in\CR$ and $\delta>0$ be such that
${\rm supp}(\varphi)\subset [a,b]\times B(\delta)$. Then, if $t\in [0,1]$, the support
of the function $(s,x)\mapsto\varphi(s-t,x)$ is contained in $[a,b+1]\times B(\delta)$.
Therefore, for any $h\in (0,1]$ we have
\begin{eqnarray*}
&\hs{5}&\int_{\CR^{d+1}}\frac{\T (h)\varphi - \varphi }{h} \, d\nu\\[1mm]
&\hs{5}=\hs{5}&\int_{[a,b+1]\times\CR^d}\frac{\T (h)\varphi - \varphi }{h} \, d\nu\\[1mm]
&\hs{5} =\hs{5} & \int_a^{b+1} \int_{\CR^d}\frac{ (G(s, s-h )\varphi (s-h ,\cdot))(x) -
(G(s, s-h) \varphi (s,\cdot ))(x)}{h} \, \mu_s(dx)\,ds \\[1mm]
&\hs{5}  & + \int_a^{b+1}\int_{\CR^d} \frac{ (G(s, s-h )\varphi (s , \cdot))(x) -
\varphi (s,x )}{h} \, \mu_s(dx)\,ds\\[1mm]
&\hs{5} =:\hs{5} & \int_a^{b+1} (I_1(s,h) + I_2(s,h)) \, ds.
\end{eqnarray*}
As far as $I_1$ is concerned, we note that
\begin{eqnarray*}
\sup_{\CR^d}\left |G(s,s-h)\left (\frac{\varphi (s-h , \cdot) - \varphi (s,\cdot)}{h}\right )\right |
\leq \sup_{\CR^{d+1} }| \varphi _t| < +\infty ,
\end{eqnarray*}
and moreover
 \[ \lim_{h\to 0^+}  \frac{\varphi (s-h , x) - \varphi (s, x)}{h}  = -  \varphi _s(s, x),
\]
the convergence being uniform in $x$. Since $f:=-
\varphi  _s \in C_c(\CR^{d+1})$,  $G(s,s-h)f$ converges uniformly to $f$ as $h\to 0^+$
by Lemma \ref{sderivative}.
Overall we see
that
\[ \lim_{h\to 0^+} I_1(s,h) = \int_{\CR^d}-  \varphi_s
(s,x) \, \mu_s(dx),
\]
and $I_1(s,h)$ is bounded by $\sup_{\CR^d}|D_t\varphi|$.

Let us consider $ I_2(s, h)$. Taking Lemma \ref{sderivative}
into account, we   write
\begin{eqnarray*}
I_2(s, h) &\hs{5} =\hs{5} & \frac{1}{h}\int_{\CR^d}  \int_{s-h}^s (G(s,r) \A(r)
\varphi (s,\cdot ))(x)dr \,\mu_s(dx) \\
&\hs{5} =\hs{5} & \frac{1}{h} \int_{s-h}^s \int_{\CR^d} (G(s,r) \A(r) \varphi
(s,\cdot))(x) \, \mu_s(dx)\,dr\\
&\hs{5} =\hs{5} & \frac{1}{h} \int_{s-h}^s \int_{\CR^d} \A(r) \varphi
(s,x) \, \mu_r(dx)\,dr,
\end{eqnarray*}
so that
\[ \lim_{h\to 0^+} I_2(s, h)  = \int_{\CR^d} \A(s) \varphi (s,x)\mu_s(dx),
\]
for almost every $s$, by the Lebesgue differentiation theorem. We
also note that
\begin{eqnarray*}
\sup_{s\in [a,b+1]}|I_2(s,h)|\le
\sup_{{r\in [a-1,b+1]}\atop{(s,x)\in {\rm supp}(\varphi)}}|(\A(r)\varphi(s,\cdot))(x)|.
\end{eqnarray*}
Hence, by the dominated convergence theorem,
\[
\lim_{h\to 0^+} \int_{\CR} (I_1(s,h) + I_2(s,h)) \, ds =
\int_{\CR}\int_{\CR^d} \left (-   \varphi _s(s,x) +
\A(s) \varphi (s,x)\right )\mu_s(dx)\,ds. \]
This proves (ii).
\end{proof}

\begin{rem}
{\rm
In view of \eqref{inf-inv} we say that $\nu$ is infinitesimally invariant, although it
is not a probability measure.
}
\end{rem}

\begin{prop}\label{strongcont}
For any $ p\in [1,+\infty)$, the semigroup $\{\T(t)\}$ extends uniquely to a
strongly continuous semigroup of positive contractions $\{\T_p(t)\}$ on
$L^p(\CR^{d+1}, \nu )$. Moreover, the infinitesimal generator of
$\{\T_p(t)\}$ is
an extension of the operator
${\mathcal G}_0:C^{1,2}_c(\CR^{d+1})\to L^p(\CR^{d+1},\nu)$
defined by
${\mathcal G}_0f={\mathcal G}f$, for any $f\in C^{1,2}_c(\CR^{d+1})$, where
${\mathcal G}$ is given by \eqref{operat-G}.
\end{prop}

\begin{proof}
Using the H\"older inequality and taking proposition \ref{transition} into account, it is immediate to check that
\begin{eqnarray*}
|(\T(t)f)(s,x)|^p
\le (\T (t) |f|^p)(s,x),\qquad\;\,(s,x)\in\CR^{d+1},\;\,t>0,
\end{eqnarray*}
for any $f \in C_c(\CR^{d+1})$. Integrating   in $\CR^{d+1}$,
we obtain
\begin{equation}
\|\T (t)f\|_{L^p(\CR^{d+1},\nu)}\le \|f\|_{L^p(\CR^{d+1},\nu)},\qquad\;\,t>0.
\label{contr}
\end{equation}
Since $C_c(\CR^{d+1})$ is dense in $L^p(\CR^{d+1},\nu)$,
estimate \eqref{contr} implies that
any operator $\T(t)$ can be extended uniquely to a bounded
operator $\T_p(t)$ which also satisfies \eqref{contr}.

Clearly, $\{\T_p(t)\}$ satisfies the semigroup law since $\{\T(t)\}$ does. It
remains to show that $\{\T_p(t)\}$ is strongly continuous.
Of course, it suffices to show that $\T_p(t)f \to f$ as
$t\to 0^+$ for all $f \in C^{1,2}_c(\CR^{d+1})$. For such $f$'s, we have
$\T_p(t)f = \T (t) f \to f$ pointwise a.e. as $t\to 0^+$ (see the proof of Lemma
\ref{behoncont}(2), where it was shown that the difference
quotients converge pointwise a.e.) and the functions $\T_p(t)f$ are uniformly bounded. The
dominated convergence theorem implies that $\{\T_p(t)\}$ is strongly continuous.

To complete the proof, let us prove that $C^{1,2}_c(\CR^{d+1})$ is contained in the
domain of the infinitesimal generator of the semigroup $\{\T_p(t)\}$. For this purpose,
we adapt the proof of Lemma \ref{behoncont}(ii). Let $a,b\in\CR$ and $\delta>0$ be such that
${\rm supp}(\varphi)\subset [a,b]\times B(\delta)$.

By Lemma \ref{sderivative} we know that
\begin{eqnarray*}
(G(s,s-t)f(s-t,\cdot))(x) - f(s-t,x) = \int_{s-t}^{s}(G(s,r)\A(r)f(s-t,\cdot))(x)\, dr,
\end{eqnarray*}
for any $(s,x)\in\CR^{d+1}$ and any $t>0$. It follows that
\begin{eqnarray*}
&\hs{5}&\left |\frac{(\T(t)f)(s,x)-f(s,x)}{t}-({\mathcal G}f)(s,x)\right |\\
&\hs{5}\le\hs{5}&\frac{1}{t}\int_{s-t}^{s}|(G(s,r)\A(r)f(s-t,\cdot))(x)-(\A(s)f(s,\cdot))(x)|\, dr\\[1mm]
&\hs{5}&+\frac{1}{t}\int_{s-t}^s| f_t(r,x)- f_t(s,x)|\,dr,
\end{eqnarray*}
for any $(s,x)\in\CR^{d+1}$ and any $t>0$.
Arguing as in the proof of Proposition \ref{prop6.1} it is immediate to check that
the function $(r,p)\mapsto (G(s,r)\A(r)f(p,\cdot))(x)$ is continuous in $\{(r,p)\in\CR^2: r\le s\}$.
Therefore,
\begin{eqnarray*}
\lim_{t\to 0^+}
(G(s,r)\A(r)f(s-t,\cdot))(x)-(\A(s)f(s,\cdot))(x)=0,
\end{eqnarray*}
for any $(s,x)\in \CR^{d+1}$.
Thus,
\begin{eqnarray*}
\lim_{t\to 0^+}\frac{(\T(t)f)(s,x)-f(s,x)}{t}=({\mathcal G}f)(s,x),\qquad\;\,(s,x)\in\CR^{d+1}.
\end{eqnarray*}
Moreover,
\begin{eqnarray*}
&\hs{5}&\sup_{(s,x)\in [a,b+1]\times B_{\delta}}\left |\frac{(\T(t)f)(s,x)-f(s,x)}{t}\right |\\[1mm]
&\hs{5}\le\hs{5}& \sup_{(r,s,x)\in [a,b+1]\times {\rm supp}(\varphi)}|(\A(r)f(s,\cdot))(x)|
+\|D_tf\|_{\infty}.
\end{eqnarray*}
Hence, the dominated convergence theorem implies that
$t^{-1}(\T(t)f-f)$ converges to ${\mathcal G}f$,  as $t\to 0^+$, in $L^p(\CR^{d+1},\nu)$ for any
$p\in [1,+\infty)$.
\end{proof}

\section{An example}

In this section we consider operators $\A(t)$ defined on smooth functions
$\varphi:\CR^d\to\CR$ by
\begin{equation}
(\A(t)\varphi)(x) = \Delta \varphi(x) + \langle b(t,x), \nabla \varphi(x) \rangle,
\label{oper-A}
\end{equation}
under the following assumptions on $b=(b_1,\ldots,b_d)$.
\begin{hyp}
\label{hyp6}
\begin{itemize}
\item[(i)]
the functions $b_j$ $(j=1,\ldots,d)$ and their first-order spatial derivatives belong to $C^{\frac{\alpha}{2},\alpha}_{\rm loc}(I\times\CR^d)$
for some $\alpha\in (0,1)$;
\item[(ii)]
the function $b(\cdot,0)$ is bounded in $I$;
\item[(iii)]
there exists a continuous function $C:I\to\CR$ such that
\begin{itemize}
\item[(a)]
$C$ is bounded from above in $I$;
\item[(b)]
$\limsup_{t\to +\infty}C(t)<0$;
\item[(c)]
$\langle \nabla_xb(t,x)\xi,\xi\rangle\le C(t)|\xi|^2,\qquad\;\,t\in I,\;\,x,\xi\in\CR^d$.
\end{itemize}
\end{itemize}
\end{hyp}

Under Hypothesis \ref{hyp6}, it is easy to check that, for any $N\in\CN$, the function $\varphi:\CR^d\to\CR$, defined by $\varphi(x)=1+|x|^{2N}$
for any $x\in\CR^d$, is a suitable Lyapunov function
for the operator $\A$ satisfying both Hypotheses \ref{hyp1}(iii) and \ref{hyp4}. Indeed,
a straightforward computation shows that
\begin{eqnarray*}
(\A(t)\varphi)(x)=|x|^{2N-2}\left (4N^2+2N(d-2)+2N\langle b(t,x),x\rangle\right ).
\end{eqnarray*}
Using Hypothesis \ref{hyp6}(iii)(c), yields
\begin{eqnarray*}
b_j(t,x)= b_j(t,0)+\int_0^1\frac{d}{ds}b_j(t,sx)ds=b_j(t,0)+\int_0^1\langle \nabla_x b_j(t,sx),x\rangle ds,
\end{eqnarray*}
so that
\begin{eqnarray*}
2\langle b(t,x),x\rangle&\hs{5}=\hs{5}&2\langle b(t,0),x\rangle +2\int_0^1\langle \nabla_xb(t,sx)x,x\rangle ds\\
&\hs{5}\le\hs{5}& 2|b(t,0)||x|+2C(t)|x|^2,
\end{eqnarray*}
for any $t\in I$ and any $x\in\CR^d$.  Hence, for any $\varepsilon>0$, we have
\begin{equation}
(\A(t)\varphi)(x)\le
(4N^2+2N(d-2))|x|^{2N-2}+2|b(t,0)||x|^{2N-1}+2NC(t)|x|^{2N}.
\label{stima-A}
\end{equation}
Since
\begin{eqnarray*}
|x|^{2N-j}\le \varepsilon |x|^{2N}+C_{\frac{2N}{j}}\varepsilon^{1-\frac{2N}{j}},\qquad\;\,x\in\CR^N,\;\,
\varepsilon>0,\;\,j=1,2,
\end{eqnarray*}
where $C_m=(m/(m-1))^{1-m}/m$,
we can rewrite \eqref{stima-A} as follows:
\begin{eqnarray}
(A(t)\varphi)(x)
&\hs{5}\le\hs{5}& \{2NC(t)+\varepsilon(2N(d-2)+4N^2)+2\varepsilon N|b(t,0)|\} |x|^{2N}\nonumber\\
&\hs{5}&+C_N\left (2N(d-2)+4N^2\right )\varepsilon^{1-N}+2NC_{2N}\varepsilon^{1-2N}|b(t,0)|
\nonumber
\\
&\hs{5} := \hs{5}&\psi_1(t) |x|^{2N} + \psi_2(t).
\label{stima-A1}
\end{eqnarray}
Hypothesis \ref{hyp1}(iii)   follows   taking $\varepsilon=1$
and, for any bounded interval $J$ compactly
supported in $I$, $\lambda_J \geq \max\{ \sup_J \psi_1, \, \sup_J \psi_2\} $.
Similarly, if we fix $\varepsilon=\varepsilon_N$ such that
\begin{eqnarray*}
\varepsilon\left (d-2+2N+|b(\cdot,0)\|_{\infty}\right )\le -\frac{1}{2}\limsup_{t\to +\infty}C(t),
\end{eqnarray*}
and $t_0\in\CR$ such that
\begin{equation}
C(t)<\frac{1}{2}\limsup_{\tau\to +\infty}C(\tau),\qquad\;\,t\ge t_0,
\label{giovedi-pioggia}
\end{equation}
then Hypothesis \ref{hyp4} is satisfied for any $t\ge t_0$.

Finally, Hypothesis \ref{hyp3}(ii) is trivially satisfied by virtue of Hypothesis \ref{hyp6}(iii)(a).
Hence, the following result follows.
\begin{thm}
Let $\A$  be defined by \eqref{oper-A} with the function $b$ satisfying Hypothesis
$\ref{hyp6}$. Then, problem
\eqref{pde} is well posed in $C_b(\CR^d)$. The corresponding evolution family $\{G(t,s)\}$
is irreducible and maps bounded measurable functions into bounded continuous functions.
Moreover, $\{G(t,s)\}$ admits an evolution system ($\mu_t)$ of measures having bounded moments
of any order $N\in\CN$ for any $t\ge t_0$, where $t_0$ is the number in \eqref{giovedi-pioggia}.
In particular, all the polynomial functions
are integrable with respect to $\mu_t$ for any $t\ge t_0$.
\end{thm}

\end{document}